\documentclass[11pt,a4paper]{article}
\usepackage{a4,amssymb,latexsym,amsmath}
\begin{document}


\newcounter{pkt}
\newenvironment{enum}{\setcounter{pkt}{0}
\begin{list}{\rm\alph{pkt})}{\usecounter{pkt}
\setlength{\topsep}{1ex}\setlength{\labelwidth}{0.5cm}
\setlength{\leftmargin}{1cm}\setlength{\labelsep}{0.25cm}
\setlength{\parsep}{-3pt}}}{\end{list}~\\[-6ex]}

\newcounter{punkt1}
\newenvironment{enum1}{\setcounter{punkt1}{0}
\begin{list}{\arabic{punkt1})}{\usecounter{punkt1}
\setlength{\topsep}{1ex}\setlength{\labelwidth}{0.6cm}
\setlength{\leftmargin}{1cm}\setlength{\labelsep}{0.25cm}
\setlength{\parsep}{-3pt}}}{\end{list}~\\[-6ex]}

\newcounter{punkt3}
\newenvironment{enum2}{\setcounter{punkt3}{0}
\begin{list}{\rm (\roman{punkt3})}{\usecounter{punkt3}
\setlength{\topsep}{1ex}\setlength{\labelwidth}{0.6cm}
\setlength{\leftmargin}{1cm}\setlength{\labelsep}{0.25cm}
\setlength{\parsep}{-3pt}}}{\end{list}~\\[-6ex]}

\newcounter{punkt2}
\newenvironment{enumbib}{\setcounter{punkt2}{0}
\begin{list}{\arabic{punkt2}.}{\usecounter{punkt2}
\setlength{\topsep}{1ex}\setlength{\labelwidth}{0.6cm}
\setlength{\leftmargin}{1cm}\setlength{\labelsep}{0.25cm}
\setlength{\parsep}{1pt}}}{\end{list}}

\newcounter{pic}
\setcounter{pic}{0}

\newcommand{\bdpm}{\begin{displaymath}}
\newcommand{\edpm}{\end{displaymath}}

\newcommand{\bea}{\begin{eqnarray}}
\newcommand{\eea}{\end{eqnarray}}

\newcommand{\beas}{\begin{eqnarray*}}
\newcommand{\eeas}{\end{eqnarray*}}

\newcommand{\ba}{\begin{array}}
\newcommand{\ea}{\end{array}}

\newenvironment{bew}{\vspace*{-0.25cm}\begin{sloppypar}\noindent{\it 
Proof.}}{\hfill\qed\end{sloppypar}\vspace*{0.15cm}}

\newtheorem{satz}{Theorem}[section]          
\newtheorem{lem}[satz]{Lemma}
\newtheorem{kor}[satz]{Corollary}
\newtheorem{prop}[satz]{Proposition}
\newtheorem{bsp}[satz]{Example}
\newtheorem{bem}[satz]{Remark}
\newtheorem{beme}[satz]{Remarks}
\newtheorem{verm}[satz]{Conjecture}

\newcommand{\brm}{\begin{rm}}
\newcommand{\erm}{\end{rm}}

\newcommand{\qed}{\hfill $\Box$}
\newcommand{\inv}{{\sf inv}}
\newcommand{\des}{{\sf des}}


\begin{center}
{\large\bf ON THE DIAGRAM OF SCHR\"ODER PERMUTATIONS}\\[1cm]
Astrid Reifegerste\\
Institut f\"ur Algebra und Geometrie\\
Otto-von-Guericke-Universit\"at Magdeburg\\
Postfach 4120, D-39016 Magdeburg, Germany\\
{\it astrid.reifegerste@mathematik.uni-magdeburg.de}\\[0.5cm]
September 17, 2002 
\end{center}
\vspace*{0.2cm}

\begin{footnotesize}
{\sc Abstract.} Egge and Mansour have recently studied permutations 
which avoid $1243$ and $2143$ regarding the occurrence of certain additional 
patterns. Some of the open questions related to their work can easily be answered 
by using permutation diagrams. Like for $132$-avoiding permutations the 
diagram approach gives insights into the structure of $\{1243,2143\}$-avoiding 
permutations that yield simple proofs for some enumerative results concerning 
forbidden patterns in such permutations.
\end{footnotesize}
\vspace*{0.5cm}


\setcounter{section}{1}\setcounter{satz}{0}

\centerline{\large{\bf 1}\hspace*{0.25cm}
{\sc Introduction}}
\vspace*{0.65cm}

Let ${\cal S}_n$ be the set of all permutations of $\{1,\ldots,n\}$. Given a 
permutation $\pi=\pi_1\cdots\pi_n\in{\cal S}_n$ and a permutation 
$\tau=\tau_1\cdots\tau_k\in{\cal S}_k$, we say that $\pi$ {\it contains the 
pattern} $\tau$ if there is a sequence $1\le i_1<i_2<\ldots<i_k\le n$ such that 
the elements $\pi_{i_1}\pi_{i_2}\cdots\pi_{i_k}$ are in the same relative order 
as $\tau_1\tau_2\cdots\tau_k$. Otherwise, $\pi$ {\it avoids the pattern} 
$\tau$, or alternatively, $\pi$ is {\it $\tau$-avoiding}. The set of $\tau$-avoiding 
permutations in ${\cal S}_n$ is denoted by ${\cal S}_n(\tau)$. 
For an arbitrary finite collection $T$ of patterns we write ${\cal S}_n(T)$ to 
denote the permutations of $\{1,\ldots,n\}$ which avoid each pattern in 
$T$.\\[2ex]
Egge and Mansour \cite{egge-mansour} studied permutations which avoid both 
$1243$ and $2143$. This work was motivated by the parallels to $132$-avoiding 
permutations. In \cite[Lem. 2 and Cor. 9]{kremer} was shown that the 
number of elements of ${\cal S}_n(1243,2143)$ is counted by the $(n-1)$st 
Schr\"oder number $r_{n-1}$. The (large) {\it Schr\"oder numbers} may be 
defined by
\bdpm
r_0:=1,\quad r_n:=r_{n-1}+\sum_{i=0}^{n-1} r_{i}r_{n-1-i}\quad\mbox{for }n\ge1.
\edpm
For this reason the authors of \cite{egge-mansour} called the permutations 
which avoid $1243$ and $2143$ {\it Schr\"oder permutations}; we will do this as 
well. (The reference to Schr\"oder numbers may be somewhat inexact because there 
are ten inequivalent pairs $(\tau_1,\tau_2)\in{\cal S}_4^2$ for which $|{\cal 
S}_n(\tau_1,\tau_2)|=r_{n-1}$, see \cite[Theo. 3]{kremer}. However, it is sufficient 
for our purposes.)
\newpage

Schr\"oder permutations are known to have a lot of properties which are analogous 
to properties of $132$-avoiding permutations. Why it needs to be so, a look at 
their diagrams shows.\\
Given a permutation $\pi\in{\cal S}_n$, we obtain its {\it diagram} $D(\pi)$ as 
follows: first let $\pi$ be represented by an $n\times n$-array with a dot in 
each of the squares $(i,\pi_i)$. Shadow all squares due south or due east of some 
dot and the dotted cell itself. The diagram $D(\pi)$ is defined as the region 
left unshaded after this procedure. A square that belongs to $D(\pi)$ we call a {\it diagram square}; a row (column) of the array that contains a 
diagram square is called a {\it diagram row} ({\it diagram column}). (The diagram is an important tool in the 
theory on the Schubert polynomial of a permutation. Schubert polynomials were 
extensively developed by Lascoux and Sch\"utzenberger. See \cite{macdonald} for a 
treatment of this work.)\\
By the construction, each of the connected components of $D(\pi)$ is a Young diagram. 
Their corners are defined to be the elements of the {\it essential set} ${\cal E}(\pi)$ of the 
permutation $\pi$. In \cite{fulton}, Fulton introduced this set which together 
with a rank function was used as a tool for algebraic treatment of Schubert 
polynomials. For any element $(i,j)\in{\cal E}(\pi)$, its {\it rank} is defined 
to be the number of dots northwest of $(i,j)$, and is denoted by $\rho(i,j)$. 
Furthermore, by ${\cal E}_r(\pi)$ we denote the set of all elements of ${\cal 
E}(\pi)$ whose rank equals $r$.\\
It is clear from the construction that the number of dots in the northwest is the same 
for all diagram squares which are connected. Consequently, we can extend the rank function on the 
diagram squares. It is a fundamental 
property of the ranked essential set of a permutation $\pi$, that it uniquely determines 
$\pi$. (This result was first proved by Fulton, see \cite[Lem. 3.10b]{fulton}; 
alternatively, an algorithm for retrieving the permutation from its 
ranked essential set was provided in \cite{eriksson-linusson}.)\\ 
Answering a question of Fulton, Eriksson and Linusson gave in 
\cite{eriksson-linusson} a characterization of all ranked sets of 
squares that arise as ranked essential set of a permutation.\\
To recover a permutation from its diagram is trivial: row by row, put a dot in 
the leftmost shaded square such that there is exactly one dot in each column.\\[2ex]
The concept should be clear from Figure 1.
 
\begin{center}
\unitlength=0.3cm
\begin{picture}(10,10)
\linethickness{0.2pt}
\multiput(0,0)(0,1){11}{\line(1,0){10}}
\multiput(0,0)(1,0){11}{\line(0,1){10}}
\put(8.5,9.5){\circle*{0.5}}
\put(3.5,8.5){\circle*{0.5}}
\put(7.5,7.5){\circle*{0.5}}
\put(9.5,6.5){\circle*{0.5}}
\put(2.5,5.5){\circle*{0.5}}
\put(0.5,4.5){\circle*{0.5}}
\put(6.5,3.5){\circle*{0.5}}
\put(5.5,2.5){\circle*{0.5}}
\put(1.5,1.5){\circle*{0.5}}
\put(4.5,0.5){\circle*{0.5}}
\multiput(8,9)(0,0.2){5}{\line(1,0){2}}
\multiput(3,8)(0,0.2){5}{\line(1,0){7}}
\multiput(3,7)(0,0.2){5}{\line(1,0){1}}\multiput(7,7)(0,0.2){5}{\line(1,0){3}}
\multiput(3,6)(0,0.2){5}{\line(1,0){1}}\multiput(7,6)(0,0.2){5}{\line(1,0){3}}
\multiput(2,5)(0,0.2){5}{\line(1,0){8}}
\multiput(0,4)(0,0.2){5}{\line(1,0){10}}
\multiput(0,3)(0,0.2){5}{\line(1,0){1}}\multiput(2,3)(0,0.2){5}{\line(1,0){2}}
\multiput(6,3)(0,0.2){5}{\line(1,0){4}}
\multiput(0,2)(0,0.2){5}{\line(1,0){1}}\multiput(2,2)(0,0.2){5}{\line(1,0){2}}
\multiput(5,2)(0,0.2){5}{\line(1,0){5}}
\multiput(0,1)(0,0.2){5}{\line(1,0){10}}
\multiput(0,0)(0,0.2){5}{\line(1,0){10}}
\linethickness{0.7pt}
\multiput(7,9)(0,1){2}{\line(1,0){1}}\multiput(7,9)(1,0){2}{\line(0,1){1}}
\multiput(2,6)(0,1){2}{\line(1,0){1}}\multiput(2,6)(1,0){2}{\line(0,1){1}}
\multiput(1,5)(0,1){2}{\line(1,0){1}}\multiput(1,5)(1,0){2}{\line(0,1){1}}
\multiput(6,6)(0,1){2}{\line(1,0){1}}\multiput(6,6)(1,0){2}{\line(0,1){1}}
\multiput(1,2)(0,1){2}{\line(1,0){1}}\multiput(1,2)(1,0){2}{\line(0,1){1}}
\multiput(4,2)(0,1){2}{\line(1,0){1}}\multiput(4,2)(1,0){2}{\line(0,1){1}}
\multiput(5,3)(0,1){2}{\line(1,0){1}}\multiput(5,3)(1,0){2}{\line(0,1){1}}
\put(1.5,5.5){\makebox(0,0)[cc]{\sf\tiny0}}
\put(2.5,6.5){\makebox(0,0)[cc]{\sf\tiny0}}
\put(7.5,9.5){\makebox(0,0)[cc]{\sf\tiny0}}
\put(6.5,6.5){\makebox(0,0)[cc]{\sf\tiny1}}
\put(1.5,2.5){\makebox(0,0)[cc]{\sf\tiny1}}
\put(4.5,2.5){\makebox(0,0)[cc]{\sf\tiny3}}
\put(5.5,3.5){\makebox(0,0)[cc]{\sf\tiny3}}
\end{picture}
\vspace*{1ex}

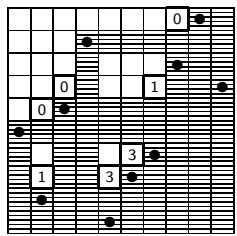
{\footnotesize{\bf Figure 1}\hspace*{0.25cm}Diagram and ranked essential set of 
$\pi=9\:4\:8\:10\:3\:1\:7\:6\:2\:5\in{\cal S}_{10}$.}
\end{center}

In \cite{reifegerste}, we used permutation diagrams to give combinatorial 
proofs for some enumerative results concerning forbidden subsequences in 
$132$-avoiding permutations. Now we develop analogues of these 
bijections. In particular, we will discuss some open problems 
which have been raised in \cite{egge-mansour}.\\[2ex]
The following section begins with a characterization of Schr\"oder permutation 
diagrams. Then we will give a surjection that takes any Schr\"oder 
permutation to a $132$-avoiding permutation of the same inversion number. On 
the other hand, a simple way to generate all Schr\"oder permutation diagrams 
from those corresponding to $132$-avoiding permutations is described.\\
Section 3 deals with additional restrictions of Schr\"oder permutations. As 
it was done for $132$-avoiding permutations we will characterize from the 
diagram the occurrence 
of increasing and decreasing subsequences of prescribed length, as well as of some 
modifications. This yields simple combinatorial proofs for 
some results appearing in \cite{egge-mansour}.\\ 
In the same reference a bijection between Schr\"oder permutations and lattice 
paths was given. Section 4 shows how the path can immediately be obtained from 
the diagram of the corresponding permutation.\\
The paper ends with some remarks about potential generalizations of its results. 

\vspace*{0.75cm}


\setcounter{section}{2}\setcounter{satz}{0}

\centerline{\large{\bf 2}\hspace*{0.25cm}
{\sc A description of Schr\"oder permutation diagrams}}
\vspace*{0.5cm}

By \cite[Theo. 2.2]{reifegerste}, $132$-avoiding permutations are precisely those 
permutations for which the diagram corresponds to a partition, or equivalently, 
for which the rank of every element of the essential set equals $0$. Analogously, we can 
characterize the elements of ${\cal S}_n(1243,2143)$.

\begin{satz} \label{characterization of Schroeder permutations}
A permutation $\pi\in{\cal S}_n$ is a Schr\"oder permutation if and only if 
every element of its essential set is of rank at most $1$.
\end{satz}

\begin{bew}
If there exists an element $(i,j)\in{\cal E}(\pi)$ (or equivalently, any 
diagram square $(i,j)$) with $\rho(i,j)\ge 2$ then, by definition, at least two dots 
appear in the northwest of $(i,j)$, say in the rows $i_1<i_2$. Obviously, the subsequence $\pi_{i_1}\pi_{i_2}\pi_i\pi_{i_3}$ is of type 
$1243$ (represented in the following figure) or $2143$ where $\pi_{i_3}=j$:
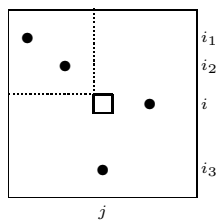
\begin{center}
\unitlength=0.25cm
\begin{picture}(10,10)
\linethickness{0.2pt}
\multiput(0,0)(0,10){2}{\line(1,0){10}}
\multiput(0,0)(10,0){2}{\line(0,1){10}}
\put(1,8.5){\circle*{0.5}}
\put(3,7){\circle*{0.5}}
\put(5,1.5){\circle*{0.5}}
\put(7.5,5){\circle*{0.5}}
\put(10.2,8.5){\makebox(0,0)[lc]{\tiny$i_1$}}
\put(10.2,7){\makebox(0,0)[lc]{\tiny$i_2$}}
\put(10.2,5){\makebox(0,0)[lc]{\tiny$i$}}
\put(10.2,1.5){\makebox(0,0)[lc]{\tiny$i_3$}}
\put(5,-0.4){\makebox(0,0)[ct]{\tiny$j$}}
\linethickness{0.7pt}
\put(4.5,5.5){\line(1,0){1}}\put(4.5,4.5){\line(1,0){1}}
\put(4.5,4.5){\line(0,1){1}}\put(5.5,4.5){\line(0,1){1}}
\linethickness{0.4pt}
\bezier{20}(4.5,5.5)(4.5,7.75)(4.5,10)
\bezier{20}(0,5.5)(2.25,5.5)(4.5,5.5)
\end{picture}
\end{center}
On the other hand, it is clear from the construction that the occurrence of a 
pattern $1243$ or $2143$ in a permutation yields a diagram corner of rank at least 
$2$.
\end{bew}
\vspace*{1ex}

We wish to describe the diagrams more precisely that arise as diagram of a Schr\"oder 
permutation. First we state two elementary properties of each permutation diagram.

\begin{lem} \label{diagram properties}
Let $\pi\in{\cal S}_n$ be an arbitrary permutation.
\begin{enum}
\item We have $i+j\le n+r$ for each $(i,j)\in{\cal E}_r(\pi)$.
\item Let $(i,j)$ be a diagram square of rank $1$ for which $(i-1,j)$ and $(i,j-1)$ do not belong to $D(\pi)$. Then 
$\pi_{i-1}=j-1$. Furthermore, for any element $(i,j)\in{\cal E}_1(\pi)$ there exists 
no square $(i',j')\in{\cal E}(\pi)$ with $i'<i$ and $j'<j$.
\end{enum}
\end{lem}
\vspace*{-0.45cm}

\begin{bew}
{\bf a)} Let $(i,j)\in{\cal E}(\pi)$ be of rank $r$. Then exactly 
$r$ indices $k<i$ satisfy $\pi_k<j$. By construction, we have $\pi_i>j$ and 
$i<\pi^{-1}_j$. Thus there exist $i-r$ integers 
$k\le i$ with $\pi_k>j$. Clearly, the number of all elements $\pi_k>j$ in $\pi$ 
equals $n-j$. This yields the restriction.\\
{\bf b)} By definition, there is exactly one dot (representing a pair $(i',j')$ 
where $\pi_{i'}=j'$) northwest of $(i,j)$. From the condition that $(i,j)$ forms the 
upper left-hand corner of a connected component of diagram squares follows 
$\pi_{i-1}<j$ and $\pi^{-1}_{j-1}<i$. Thus we have $i'=i-1$ and $j'=j-1$.\\
For the second assertion let $(i,j)\in{\cal E}_1(\pi)$. Suppose that there exists a 
diagram corner $(i',j')$ with $i'<i$ and $j'<j$. Obviously, $(i',j')$ must be 
of rank $0$, and by the first part, it is different from $(i-1,j-1)$. Thus $(i',j')$ is a corner of the Young diagram formed from the 
diagram squares that are connected with $(1,1)$. Hence $\pi_{i'+1}\le j'$ and 
$\pi^{-1}_{j'+1}\le i'$. (Note that $i'+1<i$ and $j'+1<j$; otherwise $(i,j)$ is 
not a diagram square.) Consequently, there are two dots northwest of $(i,j)$, 
contradicting to $(i,j)\in{\cal E}_1(\pi)$. 
\end{bew}

\begin{bem}
\brm
Condition a) is a part of Eriksson's and Linusson's characterization of ranked essential 
sets, see \cite[Theo. 4.1]{eriksson-linusson}.\\
In case of Schr\"oder permutations, the second claim of b) means: there are no 
diagram corners $(i,j)$ and $(i',j')$ such that $i'<i$ and $j'<j$. By \cite[Prop. 
9.6]{fulton}, just this property characterizes {\it vexillary} permutations. 
Fulton's description is an important example of characterization classes of 
permutations by the shape of their essential set. He gave a set of sufficient 
conditions that all except for one are also necessary. Later, Eriksson and 
Linusson strengthened that condition to obtain a set of necessary and 
sufficient conditions. Note that vexillary permutations can alternatively be characterized as 
$2143$-avoiding ones, see \cite[(1.27)]{macdonald}. Of course, every Schr\"oder 
permutation is vexillary.
\erm
\end{bem}

Consequently, we can answer the question when a subset of the $n^2$ squares of 
$\{1,\ldots,n\}^2$ is the essential set of a Schr\"oder permutation in ${\cal 
S}_n$. In particular, this yields a further combinatorial interpretation of 
Schr\"oder numbers.

\begin{prop} \label{fulton characterization}
For $s\ge 0$ let $i_1\ge i_2\ge\ldots\ge i_s$ and $j_1\le j_2\le\ldots\le j_s$ be positive 
integers, and let $r_1,r_2,\ldots,r_s$ be $0$ or $1$ such that
\bea
i_1-r_1>i_2-r_2>\ldots>i_s-r_s>0\quad\mbox{and}\quad 
0<j_1-r_1<j_2-r_2<\ldots<j_s-r_s.
\eea
For any $n\ge i_1+j_s$ there is a unique permutation $\pi\in{\cal S}_n$ 
with ${\cal E}(\pi)=\{(i_1,j_1),\ldots,(i_s,j_s)\}$ and $\rho(i_k,j_k)=r_k$ for 
$k=1,\ldots,s$. In particular, $\pi$ avoids $1243$ and $2143$, and every 
Schr\"oder permutation arises from a unique collection of such integers.
\end{prop}

\begin{bew}
See Fulton's \cite[Prop. 9.6]{fulton}. The condition $r_k\in\{0,1\}$ follows 
from Therorem \ref{characterization of Schroeder permutations}.
\end{bew}
\vspace*{1ex}

In \cite[Prop. 2.2]{eriksson-linusson}, the condition $n\ge i_1+j_s$ has been replaced by $i_k+j_k\le n+r_k$ for $k=1,\ldots,s$. 

\begin{kor}
\begin{enum}
\item[]
\item The $(n-1)$st Schr\"oder number $r_{n-1}$ counts the number of tripels of the integer 
sequences $i_1\ge i_2\ge\ldots\ge i_s>0$ and $0<j_1\le j_2\le\ldots\le j_s$, 
and the binary sequence $r_1,\ldots,r_s$ satisfying {\rm (1)} and $i_k+j_k\le 
n+r_k$ for all $k$.
\item The $n$th Catalan number $C_n$ counts the number of pairs of integer sequences 
$i_1>i_2>\ldots>i_s>0$ and $0<j_1<j_2<\ldots<j_s$ such that $i_k+j_k\le n$ for 
all $k$. In particular, the number of such pairs of sequences of length $s$ is 
counted by the Narayana number $N(n,s+1)$.   
\end{enum}
\end{kor}
\vspace*{-0.35cm}

\begin{bew}
The special case of $132$-avoiding permutations ($\rho(i,j)=0$ for each element $(i,j)$ of the essential 
set) in \ref{fulton characterization} yields part b). It is well known that 
$|{\cal S}_n(132)|=C_n=\frac{1}{n+1}{2n\choose n}$ for all 
$n$. The second result of b) where $N(n,s+1)=\frac{1}{n}{n\choose s}{n\choose 
s+1}$ appeared in \cite[Rem. 2.6c]{reifegerste}.
\end{bew}

Some of the results of this paper are given in terms of essential sets. Therefore we will 
describe first how one can retrieve a Schr\"oder permutation from its ranked 
essential set. In the special case of Schr\"oder permutations the retrieval 
algorithm due to Eriksson and Linusson is an evident procedure; that's why we 
will do this without any technical notation used in \cite{eriksson-linusson}.\\[2ex]
Let $\pi\in{\cal S}_n$ be a Schr\"oder permutation, and $E:={\cal E}(\pi)$ its 
essential set. Hence $E$ is a subset of labeled squares in $\{1,2,\ldots,n\}^2$ 
satisfying Proposition \ref{fulton characterization}. Let the elements of $E$ 
be represented as white labeled squares in an $n\times n$-array. (All squares that do not belong to $E$ are shaded.) 
\begin{center}
\unitlength=0.25cm
\begin{picture}(7,7)
\linethickness{0.5pt}
\multiput(0,0)(0,1){8}{\line(1,0){7}}
\multiput(0,0)(1,0){8}{\line(0,1){7}}
\put(5.5,5.5){\makebox(0,0)[cc]{\sf\tiny1}}
\put(2.5,4.5){\makebox(0,0)[cc]{\sf\tiny0}}
\put(0.5,1.5){\makebox(0,0)[cc]{\sf\tiny0}}
\put(2.5,2.5){\makebox(0,0)[cc]{\sf\tiny1}}
\linethickness{0.2pt}
\multiput(0,6)(0,0.2){5}{\line(1,0){7}}
\multiput(0,5)(0,0.2){5}{\line(1,0){5}}\multiput(6,5)(0,0.2){5}{\line(1,0){1}}
\multiput(0,4)(0,0.2){5}{\line(1,0){2}}\multiput(3,4)(0,0.2){5}{\line(1,0){4}}
\multiput(0,3)(0,0.2){5}{\line(1,0){7}}
\multiput(0,2)(0,0.2){5}{\line(1,0){2}}\multiput(3,2)(0,0.2){5}{\line(1,0){4}}
\multiput(1,1)(0,0.2){5}{\line(1,0){6}}
\multiput(0,0)(0,0.2){5}{\line(1,0){7}}
\end{picture}
\vspace*{1ex}

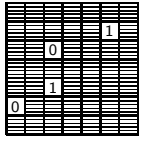
{\footnotesize{\bf Figure 2a}\hspace*{0.25cm}Ranked essential set of 
$\pi=4\:7\:5\:2\:6\:3\:1\in{\cal S}_7(1243,2143)$.}
\end{center}
(1)\hspace*{0.25cm}Colour white all squares $(i',j')$ with $i'\le i$ and $j'\le 
j$ where $(i,j)\in E$ is a square labeled with $0$. In this way we obtain {\it 
the} connected component of all diagram squares which are of rank $0$. (Note that 
the rank function can be extended on the set of diagram squares.)
\begin{center}
\unitlength=0.25cm
\begin{picture}(7,7)
\linethickness{0.5pt}
\multiput(0,0)(0,1){8}{\line(1,0){7}}
\multiput(0,0)(1,0){8}{\line(0,1){7}}
\put(5.5,5.5){\makebox(0,0)[cc]{\sf\tiny1}}
\multiput(0.5,6.5)(1,0){3}{\makebox(0,0)[cc]{\sf\tiny0}}
\multiput(0.5,5.5)(1,0){3}{\makebox(0,0)[cc]{\sf\tiny0}}
\multiput(0.5,4.5)(1,0){3}{\makebox(0,0)[cc]{\sf\tiny0}}
\put(0.5,3.5){\makebox(0,0)[cc]{\sf\tiny0}}
\put(0.5,2.5){\makebox(0,0)[cc]{\sf\tiny0}}
\put(0.5,1.5){\makebox(0,0)[cc]{\sf\tiny0}}
\put(2.5,2.5){\makebox(0,0)[cc]{\sf\tiny1}}
\linethickness{0.2pt}
\multiput(3,6)(0,0.2){5}{\line(1,0){4}}
\multiput(3,5)(0,0.2){5}{\line(1,0){2}}\multiput(6,5)(0,0.2){5}{\line(1,0){1}}
\multiput(3,4)(0,0.2){5}{\line(1,0){4}}
\multiput(1,3)(0,0.2){5}{\line(1,0){6}}
\multiput(1,2)(0,0.2){5}{\line(1,0){1}}\multiput(3,2)(0,0.2){5}{\line(1,0){4}}
\multiput(1,1)(0,0.2){5}{\line(1,0){6}}
\multiput(0,0)(0,0.2){5}{\line(1,0){7}}
\end{picture}
\vspace*{1ex}

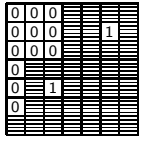
{\footnotesize{\bf Figure 2b}\hspace*{0.25cm}All diagram squares of rank $0$ are known.}
\end{center}
(2)\hspace*{0.25cm}Put a dot in each shaded square $(i,j)$ for which every square 
$(i',j')$ with $i'\le i$ and $j'\le j$, different from $(i,j)$, is a diagram square of rank 
$0$. Obviously, these dots just represent the left-to-right minima of the permutation. 
(A {\it left-to-right minimum} of a permutation $\pi$ is an element 
$\pi_i$ which is smaller than all elements to its left, i.e., $\pi_i<\pi_j$ for every $j<i$.) 
\begin{center}
\unitlength=0.25cm
\begin{picture}(7,7)
\linethickness{0.5pt}
\multiput(0,0)(0,1){8}{\line(1,0){7}}
\multiput(0,0)(1,0){8}{\line(0,1){7}}
\put(5.5,5.5){\makebox(0,0)[cc]{\sf\tiny1}}
\multiput(0.5,6.5)(1,0){3}{\makebox(0,0)[cc]{\sf\tiny0}}
\multiput(0.5,5.5)(1,0){3}{\makebox(0,0)[cc]{\sf\tiny0}}
\multiput(0.5,4.5)(1,0){3}{\makebox(0,0)[cc]{\sf\tiny0}}
\put(0.5,3.5){\makebox(0,0)[cc]{\sf\tiny0}}
\put(0.5,2.5){\makebox(0,0)[cc]{\sf\tiny0}}
\put(0.5,1.5){\makebox(0,0)[cc]{\sf\tiny0}}
\put(2.5,2.5){\makebox(0,0)[cc]{\sf\tiny1}}
\put(3.5,6.5){\circle*{0.5}}
\put(1.5,3.5){\circle*{0.5}}
\put(0.5,0.5){\circle*{0.5}}
\linethickness{0.2pt}
\multiput(3,6)(0,0.2){5}{\line(1,0){4}}
\multiput(3,5)(0,0.2){5}{\line(1,0){2}}\multiput(6,5)(0,0.2){5}{\line(1,0){1}}
\multiput(3,4)(0,0.2){5}{\line(1,0){4}}
\multiput(1,3)(0,0.2){5}{\line(1,0){6}}
\multiput(1,2)(0,0.2){5}{\line(1,0){1}}\multiput(3,2)(0,0.2){5}{\line(1,0){4}}
\multiput(1,1)(0,0.2){5}{\line(1,0){6}}
\multiput(0,0)(0,0.2){5}{\line(1,0){7}}
\end{picture}
\vspace*{1ex}

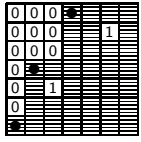
{\footnotesize{\bf Figure 2c}\hspace*{0.25cm}All dotted squares connected with 
a diagram square of rank $0$ are known.}
\end{center}
(3)\hspace*{0.25cm}For each dot contained in a square $(i,j)$, colour white 
all squares $(i'',j'')$ with $i<i''\le i'$ and $j<j''\le j'$ where $(i',j')\in 
E$ is a square labeled with $1$. By this step, all diagram squares of rank $1$ are 
obtained. (Note that all squares which are situated in the southeast area of a 
given dot belong to the same connected component.)
\begin{center}
\unitlength=0.25cm
\begin{picture}(7,7)
\linethickness{0.5pt}
\multiput(0,0)(0,1){8}{\line(1,0){7}}
\multiput(0,0)(1,0){8}{\line(0,1){7}}
\multiput(4.5,5.5)(1,0){2}{\makebox(0,0)[cc]{\sf\tiny1}}
\multiput(0.5,6.5)(1,0){3}{\makebox(0,0)[cc]{\sf\tiny0}}
\multiput(0.5,5.5)(1,0){3}{\makebox(0,0)[cc]{\sf\tiny0}}
\multiput(0.5,4.5)(1,0){3}{\makebox(0,0)[cc]{\sf\tiny0}}
\put(0.5,3.5){\makebox(0,0)[cc]{\sf\tiny0}}
\put(0.5,2.5){\makebox(0,0)[cc]{\sf\tiny0}}
\put(0.5,1.5){\makebox(0,0)[cc]{\sf\tiny0}}
\put(2.5,2.5){\makebox(0,0)[cc]{\sf\tiny1}}
\put(3.5,6.5){\circle*{0.5}}
\put(1.5,3.5){\circle*{0.5}}
\put(0.5,0.5){\circle*{0.5}}
\linethickness{0.2pt}
\multiput(3,6)(0,0.2){5}{\line(1,0){4}}
\multiput(3,5)(0,0.2){5}{\line(1,0){1}}\multiput(6,5)(0,0.2){5}{\line(1,0){1}}
\multiput(3,4)(0,0.2){5}{\line(1,0){4}}
\multiput(1,3)(0,0.2){5}{\line(1,0){6}}
\multiput(1,2)(0,0.2){5}{\line(1,0){1}}\multiput(3,2)(0,0.2){5}{\line(1,0){4}}
\multiput(1,1)(0,0.2){5}{\line(1,0){6}}
\multiput(0,0)(0,0.2){5}{\line(1,0){7}}
\end{picture}
\vspace*{1ex}

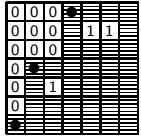
{\footnotesize{\bf Figure 2d}\hspace*{0.25cm}The diagram is completed.}
\end{center}
(4)\hspace*{0.25cm}Row by row, if no dot exists in the row, put such a one 
in the leftmost shaded square such that there is exactly one dot in each 
column. Now the permutation can read off from the array.
\begin{center}
\unitlength=0.25cm
\begin{picture}(7,7)
\linethickness{0.5pt}
\multiput(0,0)(0,1){8}{\line(1,0){7}}
\multiput(0,0)(1,0){8}{\line(0,1){7}}
\multiput(4.5,5.5)(1,0){2}{\makebox(0,0)[cc]{\sf\tiny1}}
\multiput(0.5,6.5)(1,0){3}{\makebox(0,0)[cc]{\sf\tiny0}}
\multiput(0.5,5.5)(1,0){3}{\makebox(0,0)[cc]{\sf\tiny0}}
\multiput(0.5,4.5)(1,0){3}{\makebox(0,0)[cc]{\sf\tiny0}}
\put(0.5,3.5){\makebox(0,0)[cc]{\sf\tiny0}}
\put(0.5,2.5){\makebox(0,0)[cc]{\sf\tiny0}}
\put(0.5,1.5){\makebox(0,0)[cc]{\sf\tiny0}}
\put(2.5,2.5){\makebox(0,0)[cc]{\sf\tiny1}}
\put(3.5,6.5){\circle*{0.5}}
\put(6.5,5.5){\circle*{0.5}}
\put(4.5,4.5){\circle*{0.5}}                     
\put(1.5,3.5){\circle*{0.5}}
\put(5.5,2.5){\circle*{0.5}}
\put(2.5,1.5){\circle*{0.5}}
\put(0.5,0.5){\circle*{0.5}}
\linethickness{0.2pt}
\multiput(3,6)(0,0.2){5}{\line(1,0){4}}
\multiput(3,5)(0,0.2){5}{\line(1,0){1}}\multiput(6,5)(0,0.2){5}{\line(1,0){1}}
\multiput(3,4)(0,0.2){5}{\line(1,0){4}}
\multiput(1,3)(0,0.2){5}{\line(1,0){6}}
\multiput(1,2)(0,0.2){5}{\line(1,0){1}}\multiput(3,2)(0,0.2){5}{\line(1,0){4}}
\multiput(1,1)(0,0.2){5}{\line(1,0){6}}
\multiput(0,0)(0,0.2){5}{\line(1,0){7}}
\end{picture}
\vspace*{1ex}

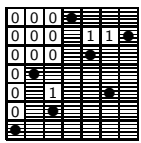
{\footnotesize{\bf Figure 2e}\hspace*{0.25cm}The permutation 
$\pi=4\:7\:5\:2\:6\:3\:1$ is recovered.}
\end{center}

The following transformation explains the close connection between $132$-avoiding permutations
and Schr\"oder permutations. 

\begin{prop} \label{transformation}
Let $\pi\in{\cal S}_n$ be a Schr\"oder permutation.
Let ${\cal E}^*(\pi)$ be the set which we obtain from ${\cal E}(\pi)$ by replacing 
every element $(i,j)\in{\cal E}_1(\pi)$ by $(i-1,j-1)$ and defining it to be of rank $0$. Then 
${\cal E}^*(\pi)$ is an essential set. (In particular, ${\cal E}^*(\pi)$ is the essential set of a 
$132$-avoiding permutation.)  
\end{prop}

\begin{bew}
Let ${\cal E}(\pi)=\{(i_1,j_1),\ldots,(i_s,j_s)\}$. We may assume that 
$\rho(i_k,j_k)=1$ for any $k$, otherwise the assertion is trivial. Set 
$i'_k:=i_k-1,\;j'_k:=j_k-1,\;r'_k:=r_k-1=0$, and check Proposition \ref{fulton 
characterization} for $E={\cal 
E}(\pi)\cup\{(i'_k,j'_k)\}\setminus\{(i_k,j_k)\}$. Evidently, all the 
conditions are satisfied (we have $i_k-r_k=i'_k-r'_k,\;j_k-r_k=j'_k-r'_k$).
\end{bew}

\begin{bsp}
\brm
Let $\pi=4\:7\:5\:2\:6\:3\:1\in{\cal S}_7(1243,2143)$. Then 
the transformation ${\cal E}(\pi)\mapsto{\cal E}^*(\pi)$ yields the essential set of 
$\sigma=6\:4\:5\:3\:2\:7\:1\in{\cal S}_7(132)$:
\begin{center}
\unitlength=0.25cm
\begin{picture}(20,7)
\linethickness{0.3pt}
\multiput(0,0)(0,1){8}{\line(1,0){7}}
\multiput(0,0)(1,0){8}{\line(0,1){7}}
\put(5.5,5.5){\makebox(0,0)[cc]{\sf\tiny1}}
\put(2.5,4.5){\makebox(0,0)[cc]{\sf\tiny0}}
\put(0.5,1.5){\makebox(0,0)[cc]{\sf\tiny0}}
\put(2.5,2.5){\makebox(0,0)[cc]{\sf\tiny1}}
\multiput(13,0)(0,1){8}{\line(1,0){7}}
\multiput(13,0)(1,0){8}{\line(0,1){7}}
\put(17.5,6.5){\makebox(0,0)[cc]{\sf\tiny0}}
\put(15.5,4.5){\makebox(0,0)[cc]{\sf\tiny0}}
\put(13.5,1.5){\makebox(0,0)[cc]{\sf\tiny0}}
\put(14.5,3.5){\makebox(0,0)[cc]{\sf\tiny0}}
\linethickness{0.1pt}
\multiput(3,6)(0,0.2){5}{\line(1,0){4}}
\multiput(3,5)(0,0.2){5}{\line(1,0){1}}\multiput(6,5)(0,0.2){5}{\line(1,0){1}}
\multiput(3,4)(0,0.2){5}{\line(1,0){4}}
\multiput(1,3)(0,0.2){5}{\line(1,0){6}}
\multiput(1,2)(0,0.2){5}{\line(1,0){1}}\multiput(3,2)(0,0.2){5}{\line(1,0){4}}
\multiput(1,1)(0,0.2){5}{\line(1,0){6}}
\multiput(0,0)(0,0.2){5}{\line(1,0){7}}
\multiput(18,6)(0,0.2){5}{\line(1,0){2}}
\multiput(16,5)(0,0.2){5}{\line(1,0){4}}
\multiput(16,4)(0,0.2){5}{\line(1,0){4}}
\multiput(15,3)(0,0.2){5}{\line(1,0){5}}
\multiput(14,2)(0,0.2){5}{\line(1,0){6}}
\multiput(14,1)(0,0.2){5}{\line(1,0){6}}
\multiput(13,0)(0,0.2){5}{\line(1,0){7}}
\linethickness{0.7pt}
\multiput(5,5)(0,1){2}{\line(1,0){1}}\multiput(5,5)(1,0){2}{\line(0,1){1}}
\multiput(2,4)(0,1){2}{\line(1,0){1}}\multiput(2,4)(1,0){2}{\line(0,1){1}}
\multiput(0,1)(0,1){2}{\line(1,0){1}}\multiput(0,1)(1,0){2}{\line(0,1){1}}
\multiput(2,2)(0,1){2}{\line(1,0){1}}\multiput(2,2)(1,0){2}{\line(0,1){1}}
\multiput(17,6)(0,1){2}{\line(1,0){1}}\multiput(17,6)(1,0){2}{\line(0,1){1}}
\multiput(15,4)(0,1){2}{\line(1,0){1}}\multiput(15,4)(1,0){2}{\line(0,1){1}}
\multiput(13,1)(0,1){2}{\line(1,0){1}}\multiput(13,1)(1,0){2}{\line(0,1){1}}
\multiput(14,3)(0,1){2}{\line(1,0){1}}\multiput(14,3)(1,0){2}{\line(0,1){1}}
\put(10,3.5){\makebox(0,0)[cc]{$\longrightarrow$}}
\end{picture}
\vspace*{1ex}

{\footnotesize{\bf Figure 3}\hspace*{0.25cm}On the left the diagram of $\pi$; 
on the right the diagram of $\sigma$.}
\end{center}
\erm
\end{bsp}
\vspace*{0.35cm}                                                               

Let $\phi:{\cal S}_n(1243,2143)\to{\cal S}_n$ be the map which takes any 
Schr\"oder permutation $\pi$ to the permutation whose essential set equals 
${\cal E}^*(\pi)$. Obviously, $\phi$ is a surjection to ${\cal S}_n(132)$.\\ 
It follows from Lemma \ref{diagram properties}b and the retrieval procedure 
that $D(\pi)$ and $D(\phi(\pi))$ have the same number of squares. By \cite[(1.21)]{macdonald}, for any permutation $\pi\in{\cal S}_n$ the number of diagram squares is equal 
to the number of inversions $\inv(\pi)$ of $\pi$. Thus we have $\inv(\pi)=\inv(\phi(\pi))$ for 
every $\pi\in{\cal S}_n(1243,2143)$. Furthermore, Fulton observed in \cite{fulton} that a 
permutation $\pi\in{\cal S}_n$ has a descent at position $i$ if and only if 
there exists a diagram corner in the $i$th row of the $n\times n$-array 
representing $\pi$. 
(An integer $i\in\{1,\ldots,n-1\}$ for which $\pi_i>\pi_{i+1}$ is called a 
{\it descent} of $\pi\in{\cal S}_n$. The number of descents of $\pi$ is denoted 
by $\des(\pi)$.) Lemma \ref{diagram properties}b implies that $\des(\pi)\le\des(\phi(\pi))$ for 
each $\pi\in{\cal S}_n(1243,2143)$.\\
The left-to-right minima of a permutation $\pi\in{\cal S}_n$ are represented by 
such dots $(i,j)$ for which $(i-1,j)$ or $(i,j-1)$ are diagram squares of rank 
$0$. To master the case $\pi_1=1$ we assume that $(0,1)$ is a diagram 
square of rank $0$. Consequently, every left-to-right minimum of $\pi\in{\cal 
S}_n(1243,2143)$ is also such a one for $\phi(\pi)$.\\[2ex]
In \cite[Theo. 5.1]{reifegerste} we have shown that the number of subsequences of 
type $132$ in an arbitrary permutation is equal to the sum of ranks of all 
diagram squares. For Schr\"oder permutations this value is just 
the number of all diagram squares of rank $1$.\\
The conversion of the above transformation is a simple way to construct 
Schr\"oder permutations which contain a prescribed number of occurrences of the pattern $132$.\\ 
Given the essential set of a $132$-avoiding permutation $\sigma\in{\cal S}_n(132)$ 
(recall that ${\cal E}(\sigma)$ is the corner set of a Young diagram fitting in 
$(n-1,n-2,\ldots,1)$; all elements are of rank $0$), we replace some elements $(i,j)\in{\cal 
E}(\sigma)$ by $(i+1,j+1)$ and increase their label by $1$. It follows from 
Proposition \ref{fulton characterization} and Lemma \ref{diagram properties}a 
that the resulting set is an essential set of a Schr\"oder permutation in 
${\cal S}_n$ if and only if we have $i+j<n$ for all replaced elements $(i,j)$.\\  
For instance, from $\sigma=6\:4\:5\:3\:2\:7\:1\in{\cal S}_7(132)$ we obtain:
\begin{center}
\unitlength=0.2cm
\begin{picture}(77,8)
\linethickness{0.3pt}
\multiput(0,1)(0,1){8}{\line(1,0){7}}
\multiput(0,1)(1,0){8}{\line(0,1){7}}
\multiput(10,1)(0,1){8}{\line(1,0){7}}
\multiput(10,1)(1,0){8}{\line(0,1){7}}
\multiput(20,1)(0,1){8}{\line(1,0){7}}
\multiput(20,1)(1,0){8}{\line(0,1){7}}
\multiput(30,1)(0,1){8}{\line(1,0){7}}
\multiput(30,1)(1,0){8}{\line(0,1){7}}
\multiput(40,1)(0,1){8}{\line(1,0){7}}
\multiput(40,1)(1,0){8}{\line(0,1){7}}
\multiput(50,1)(0,1){8}{\line(1,0){7}}
\multiput(50,1)(1,0){8}{\line(0,1){7}}
\multiput(60,1)(0,1){8}{\line(1,0){7}}
\multiput(60,1)(1,0){8}{\line(0,1){7}}
\multiput(70,1)(0,1){8}{\line(1,0){7}}
\multiput(70,1)(1,0){8}{\line(0,1){7}}
\put(4.5,7.5){\makebox(0,0)[cc]{\sf\tiny0}}
\put(2.5,5.5){\makebox(0,0)[cc]{\sf\tiny0}}
\put(1.5,4.5){\makebox(0,0)[cc]{\sf\tiny0}}
\put(0.5,2.5){\makebox(0,0)[cc]{\sf\tiny0}}
\put(15.5,6.5){\makebox(0,0)[cc]{\sf\tiny1}}
\put(12.5,5.5){\makebox(0,0)[cc]{\sf\tiny0}}
\put(11.5,4.5){\makebox(0,0)[cc]{\sf\tiny0}}
\put(10.5,2.5){\makebox(0,0)[cc]{\sf\tiny0}}
\put(24.5,7.5){\makebox(0,0)[cc]{\sf\tiny0}}
\put(23.5,4.5){\makebox(0,0)[cc]{\sf\tiny1}}
\put(21.5,4.5){\makebox(0,0)[cc]{\sf\tiny0}}
\put(20.5,2.5){\makebox(0,0)[cc]{\sf\tiny0}}
\put(34.5,7.5){\makebox(0,0)[cc]{\sf\tiny0}}
\put(32.5,5.5){\makebox(0,0)[cc]{\sf\tiny0}}
\put(32.5,3.5){\makebox(0,0)[cc]{\sf\tiny1}}
\put(30.5,2.5){\makebox(0,0)[cc]{\sf\tiny0}}
\put(45.5,6.5){\makebox(0,0)[cc]{\sf\tiny1}}
\put(43.5,4.5){\makebox(0,0)[cc]{\sf\tiny1}}
\put(41.5,4.5){\makebox(0,0)[cc]{\sf\tiny0}}
\put(40.5,2.5){\makebox(0,0)[cc]{\sf\tiny0}}
\put(55.5,6.5){\makebox(0,0)[cc]{\sf\tiny1}}
\put(52.5,5.5){\makebox(0,0)[cc]{\sf\tiny0}}
\put(52.5,3.5){\makebox(0,0)[cc]{\sf\tiny1}}
\put(50.5,2.5){\makebox(0,0)[cc]{\sf\tiny0}}
\put(64.5,7.5){\makebox(0,0)[cc]{\sf\tiny0}}
\put(63.5,4.5){\makebox(0,0)[cc]{\sf\tiny1}}
\put(62.5,3.5){\makebox(0,0)[cc]{\sf\tiny1}}
\put(60.5,2.5){\makebox(0,0)[cc]{\sf\tiny0}}
\put(75.5,6.5){\makebox(0,0)[cc]{\sf\tiny1}}
\put(73.5,4.5){\makebox(0,0)[cc]{\sf\tiny1}}
\put(72.5,3.5){\makebox(0,0)[cc]{\sf\tiny1}}
\put(70.5,2.5){\makebox(0,0)[cc]{\sf\tiny0}}
\linethickness{0.1pt}
\multiput(5,7)(0,0.2){5}{\line(1,0){2}}\multiput(3,6)(0,0.2){5}{\line(1,0){4}}
\multiput(3,5)(0,0.2){5}{\line(1,0){4}}\multiput(2,4)(0,0.2){5}{\line(1,0){5}}
\multiput(1,3)(0,0.2){5}{\line(1,0){6}}\multiput(1,2)(0,0.2){5}{\line(1,0){6}}
\multiput(0,1)(0,0.2){5}{\line(1,0){7}}
\multiput(13,7)(0,0.2){5}{\line(1,0){4}}\multiput(13,6)(0,0.2){5}{\line(1,0){1}}
\multiput(16,6)(0,0.2){5}{\line(1,0){1}}\multiput(13,5)(0,0.2){5}{\line(1,0){4}}
\multiput(12,4)(0,0.2){5}{\line(1,0){5}}\multiput(11,3)(0,0.2){5}{\line(1,0){6}}
\multiput(11,2)(0,0.2){5}{\line(1,0){6}}\multiput(10,1)(0,0.2){5}{\line(1,0){7}}
\multiput(25,7)(0,0.2){5}{\line(1,0){2}}\multiput(22,6)(0,0.2){5}{\line(1,0){5}}
\multiput(22,5)(0,0.2){5}{\line(1,0){1}}\multiput(24,5)(0,0.2){5}{\line(1,0){3}}
\multiput(22,4)(0,0.2){5}{\line(1,0){1}}\multiput(24,4)(0,0.2){5}{\line(1,0){3}}
\multiput(21,3)(0,0.2){5}{\line(1,0){6}}\multiput(21,2)(0,0.2){5}{\line(1,0){6}}
\multiput(20,1)(0,0.2){5}{\line(1,0){7}}
\multiput(35,7)(0,0.2){5}{\line(1,0){2}}\multiput(33,6)(0,0.2){5}{\line(1,0){4}}
\multiput(33,5)(0,0.2){5}{\line(1,0){4}}\multiput(31,4)(0,0.2){5}{\line(1,0){6}}
\multiput(31,3)(0,0.2){5}{\line(1,0){1}}\multiput(33,3)(0,0.2){5}{\line(1,0){4}}
\multiput(31,2)(0,0.2){5}{\line(1,0){6}}\multiput(30,1)(0,0.2){5}{\line(1,0){7}}
\multiput(42,7)(0,0.2){5}{\line(1,0){5}}\multiput(42,6)(0,0.2){5}{\line(1,0){1}}
\multiput(46,6)(0,0.2){5}{\line(1,0){1}}\multiput(42,5)(0,0.2){5}{\line(1,0){1}}
\multiput(44,5)(0,0.2){5}{\line(1,0){3}}\multiput(42,4)(0,0.2){5}{\line(1,0){1}}
\multiput(44,4)(0,0.2){5}{\line(1,0){3}}\multiput(41,3)(0,0.2){5}{\line(1,0){6}}
\multiput(41,2)(0,0.2){5}{\line(1,0){6}}\multiput(40,1)(0,0.2){5}{\line(1,0){7}}
\multiput(53,7)(0,0.2){5}{\line(1,0){4}}\multiput(53,6)(0,0.2){5}{\line(1,0){1}}
\multiput(56,6)(0,0.2){5}{\line(1,0){1}}\multiput(53,5)(0,0.2){5}{\line(1,0){4}}
\multiput(51,4)(0,0.2){5}{\line(1,0){6}}\multiput(51,3)(0,0.2){5}{\line(1,0){1}}
\multiput(53,3)(0,0.2){5}{\line(1,0){4}}\multiput(51,2)(0,0.2){5}{\line(1,0){6}}
\multiput(50,1)(0,0.2){5}{\line(1,0){7}}
\multiput(65,7)(0,0.2){5}{\line(1,0){2}}\multiput(61,6)(0,0.2){5}{\line(1,0){6}}
\multiput(61,5)(0,0.2){5}{\line(1,0){1}}\multiput(64,5)(0,0.2){5}{\line(1,0){3}}
\multiput(61,4)(0,0.2){5}{\line(1,0){1}}\multiput(64,4)(0,0.2){5}{\line(1,0){3}}
\multiput(61,3)(0,0.2){5}{\line(1,0){1}}\multiput(63,3)(0,0.2){5}{\line(1,0){4}}
\multiput(61,2)(0,0.2){5}{\line(1,0){6}}\multiput(60,1)(0,0.2){5}{\line(1,0){7}}
\multiput(71,7)(0,0.2){5}{\line(1,0){6}}\multiput(71,6)(0,0.2){5}{\line(1,0){1}}
\multiput(76,6)(0,0.2){5}{\line(1,0){1}}\multiput(71,5)(0,0.2){5}{\line(1,0){1}}
\multiput(74,5)(0,0.2){5}{\line(1,0){3}}\multiput(71,4)(0,0.2){5}{\line(1,0){1}}
\multiput(74,4)(0,0.2){5}{\line(1,0){3}}\multiput(71,3)(0,0.2){5}{\line(1,0){1}}
\multiput(73,3)(0,0.2){5}{\line(1,0){4}}\multiput(71,2)(0,0.2){5}{\line(1,0){6}}
\multiput(70,1)(0,0.2){5}{\line(1,0){7}}
\put(3.5,0){\makebox(0,0)[cc]{\tiny$\pi_1=\sigma$}}
\put(13.5,0){\makebox(0,0)[cc]{\tiny$\pi_2=4753261$}}
\put(23.5,0){\makebox(0,0)[cc]{\tiny$\pi_3=6357241$}}
\put(33.5,0){\makebox(0,0)[cc]{\tiny$\pi_4=6452731$}}
\put(43.5,0){\makebox(0,0)[cc]{\tiny$\pi_5=3756241$}}
\put(53.5,0){\makebox(0,0)[cc]{\tiny$\pi_6=4752631$}}
\put(63.5,0){\makebox(0,0)[cc]{\tiny$\pi_7=6257431$}}
\put(73.5,0){\makebox(0,0)[cc]{\tiny$\pi_8=2756431$}}
\end{picture}
\vspace*{1ex}

{\footnotesize{\bf Figure 4}\hspace*{0.25cm}(All the) Schr\"oder permutations obtained 
from $\sigma$.}
\end{center}

Obviously, these are all the Schr\"oder permutations which can 
be constructed in this way, that is, whose image with respect to $\phi$ equals 
$\sigma$. Note that $\inv(\sigma)=15=\inv(\pi_i)$ for $i=1,\ldots,8$. 

\begin{prop}
Let $\sigma\in{\cal S}_n(132)$, and let $s$ be the number of elements $(i,j)\in{\cal E}(\sigma)$ 
satisfying $i+j<n$. Then there exist $2^s$ Schr\"oder permutations $\pi\in{\cal S}_n$ for which 
$\phi(\pi)=\sigma$.
\end{prop}

\begin{bew}
This follows from the preceding discussion.
\end{bew}
\vspace*{1ex}

In \cite[Cor. 3.7]{reifegerste}, we have enumerated the Young diagrams fitting in 
$(n-1,n-2,\ldots,1)$ according to the number of their corners in the diagonal $i+j=n$. 
The number of such diagrams with $k\ge0$ corners $(i,n-i)$ equals the ballot number
$b(n-1,n-1-k)=\frac{k+1}{2n-1-k}{2n-1-k\choose n}$.\\
Now we are interested in the distribution of corners outside that diagonal.

\begin{prop}
Let $c(n-1,k)$ be the number of Young diagrams fitting in the shape 
$(n-1,n-2,\ldots,1)$ with $k\ge1$ corners satisfying $i+j<n$. Then we have
\bdpm
c(n-1,k)=\sum_{i=1}^{n-1-k} \frac{i}{n-i}{n-i\choose k}{n-1\choose k+i}.
\edpm
Furthermore, there are $2^{n-1}$ such diagrams with no corner outside the 
diagonal $i+j=n$.
\end{prop}

\begin{bew}
Consider the Young diagram as being contained in an $n\times n$-rectangle, and 
consider the lattice path from the upper right-hand to the lower left-hand 
corners of the rectangle that travels along the boundary of the diagram. 
Defining the rectangle diagonal to be the $x$-axis with origin in the lower 
left-hand corner, we obtain a Dyck path of length $2n$, that is, a lattice path from $(0,0)$ to $(2n,0)$ which never falls below 
the $x$-axis. (In \cite{reifegerste}, we have noted that the lattice path 
resulting from the diagram of a $132$-avoiding permutation 
$\pi\in{\cal S}_n$ in this way is just the Dyck path corresponding 
to $\pi$ according a bijection proposed by Krattenthaler in 
\cite{krattenthaler}.) In terms of Dyck paths, a diagram corner satisfying 
$i+j<n$ means a valley at a level greater than $0$ (where the $x$-axis marks 
the $0$-level). The distribution of the number of these valleys was given in 
\cite[Sect. 6.11]{deutsch}.
\end{bew}
\vspace*{1ex}

The previous both propositions immediately yield an explicit description for the Schr\"oder numbers.

\begin{kor}
For $n\ge 0$ we have $r_n=2^n+\sum_{k=1}^{n-1} 2^k c(n,k)$.
\end{kor}

\begin{bem}
\brm
Another one is $r_n=\sum_{k=0}^n {2n-k\choose k}C_{n-k}$ where 
$C_n=\frac{1}{n+1}{2n\choose n}$ denotes the $n$th Catalan number. This 
formula follows directly from an interpretation in terms of lattice paths, see 
\cite[Exc. 6.19 and 6.39]{stanley}.
\erm
\end{bem}
\vspace*{0.75cm}


\setcounter{section}{3}\setcounter{satz}{0}

\centerline{\large{\bf 3}\hspace*{0.25cm}
{\sc Forbidden subsequences in Schr\"oder permutations}}
\vspace*{0.5cm}
            
In this section we will demonstrate that diagrams can be used to obtain simple 
proofs for enumerative results concerning certain restrictions of Schr\"oder 
permutations. Most of numbers $|{\cal S}_n(1243,2143,\tau)|$ appearing below 
are known from their analytical derivation in \cite{egge-mansour}.\\[2ex]
For the following investigation, only one case is really of interest: the 
essential set of $\pi\in{\cal S}_n(1243,2143)$ contains both elements of rank 
$0$ and $1$. If ${\cal E}_1(\pi)=\emptyset$ then $\pi$ avoids $132$, and all 
has been done in \cite{reifegerste}. If there is no diagram corner of rank $0$ then we have 
$\pi_1=1$, and $\pi_2\cdots\pi_n$ can be identified with a permutation in ${\cal 
S}_{n-1}(132)$. In particular, these permutations contain as many 
subsequences of type $21$ (inversions) as of type $132$. (Note that the number 
of the first ones equals the number of all diagram squares, and the number of the 
latter counts all diagram squares of rank $1$).\\[2ex]
Let us start with the consideration of increasing subsequences. In \cite[Theo. 4.1b]{reifegerste} we proved that a permutation 
$\pi\in{\cal S}_n(132)$ avoids the pattern $12\cdots k$ if and only if its 
diagram contains $(n+1-k,n-k,\ldots,1)$. (Recall that in case of $132$-avoiding 
permutations the diagram corresponds to a Young diagram fitting in the shape 
$(n-1,n-2,\ldots,1)$.) This condition will be useful for Schr\"oder 
permutations as well. 

\begin{satz} \label{12...k-avoiding}
Let $\pi\in{\cal S}_n(1243,2143)$ be a Schr\"oder permutation. 
Then $\pi$ avoids $12\cdots k$ for any $k\ge 1$ if and only if 
$\phi(\pi)$ avoids $12\cdots k$. 
\end{satz}

\begin{bew}
We may assume that the essential set ${\cal E}(\pi)$ contains at least one 
element, say $(i,j)$, of rank $1$; otherwise the assertion is trivial. The 
proof of \ref{transformation} implies that the set ${\cal E}'(\pi):={\cal 
E}(\pi)\cup\{(i-1,j-1)\}\setminus\{(i,j)\}$ is the essential set of a 
Schr\"oder permutation again. The rank of $(i-1,j-1)$ is defined as $0$. 
(Successive determining yields the set ${\cal E}^*(\pi)$ stated in 
Proposition \ref{transformation}.) Now we consider which consequences for the corresponding permutation result 
from this transformation.\\
Let $\sigma\in{\cal S}_n(1243,2143)$ such that ${\cal E}(\sigma)= {\cal 
E}'(\pi)$. Then $\sigma$ differs from $\pi$ at exactly three positions. Let 
$\pi_{i_1}$ be the element represented by the only dot in the northwest of 
$(i,j)$. Furthermore let $\pi_{i_2}=j$. Then we have 
$\sigma_i=\pi_{i_1},\;\sigma_{i_2}=\pi_i,\;\sigma_{i_1}=\pi_{i_2}$, and 
$\sigma_k=\pi_k$ for all $k$, different from $i,i_1,i_2$. The proof for this 
fact, we find in the retrieval procedure given in Section 2. (For a better understanding 
it is helpful to consider simultaneously the example following the proof.)
\begin{enum1}
\item The element $\pi_{i_1}$ is a left-to-right minimum of $\pi$. All the squares due 
north or due west of the dot $(i_1,\pi_{i_1})$ are diagram squares of rank $0$. 
Let $i'$ be the smallest integer greater than $i_1$ such that a corner of rank 
$0$ appears in row $i'$. (If such a corner does not exist, set $i'=\infty$.) By 
Lemma \ref{diagram properties}b, we have $i\le i'$. In the array representing 
$\sigma$, the square $(i-1,j-1)$ forms a corner of the $0$-component, and the 
next one appears in row $i'$. Thus the dot representing $\sigma_i$ is contained 
in column $\pi_{i_1}$.
\item By the transformation, all squares $(i'',j'')$ for which $i_1<i''\le 
i$ and $\pi_{i_1}<j''\le j$ are moved northwestwards. Let $j'$ be the column 
index of the corner of rank $0$ which appears in row $i_1-1$. (Note that for 
$i_1>1$ such a one has to exist since $\pi_{i_1}$ is a left-to-right minimum. 
For $i_1=1$ set $j'=\infty$.) 
By Lemma \ref{diagram properties}b again, we have $j\le j'$. Hence in the 
array of $\sigma$ the square $(i_1,j)$ is dotted.
\item Now all dots $(i',\sigma_{i'})$ with $i'\le i$ are fixed. It follows from 
the construction that these dots are just $(i',\pi_{i'})$ if $i'\not=i_1,i$. Since all diagram 
squares south of row $i$ appear at the same position in $D(\sigma)$, the only 
possible position for the missing dot in row $i_2$ is $(i_2,\pi_i)$. For all 
other indices $k$ we have $\sigma_k=\pi_k$. (Note that $\pi_i>\pi_{i_1}$ and 
$\pi_i>\pi_{i_2}$.)
\end{enum1}
Consequently, if $\pi$ contains any increasing subsequence of length $k$, the 
permutation $\sigma$ contains such a sequence as well, and vice versa: by 
definition and step 2), the elements $\pi_{i_1}$ and $\sigma_{i_1}$ are 
left-to-right minima of $\pi$ and $\sigma$, respectively. If these elements 
occur in an increasing subsequence then as the first term. Obviously, all 
the elements $\pi_{i_1+1},\ldots,\pi_{i-1}$ are greater than $\pi_{i_2}$ 
($=\sigma_{i_1}$). Furthermore we have $\pi_k<\pi_{i_2}$ for 
$k=i+1,\ldots,i_2-1$. (Note that there exists no diagram square in the 
southeast area of $(i,j)$.) Thus, and since $\pi_{i_1}<\pi_{i_2}<\pi_i$ each 
increasing subsequence in $\pi_{i_1}\pi_{i_1+1}\cdots\pi_{i_2}$ corresponds to 
an increasing subsequence of the same length in 
$\sigma_{i_1}\sigma_{i_1+1}\cdots\sigma_{i_2}$.\\
Using the arguments successively (until the permutation $\phi(\pi)$ is obtained) proves the assertion of the 
theorem.
\end{bew}
\vspace*{1ex}

\begin{bsp}
\brm
For $\pi=5\:9\:8\:10\:4\:2\:6\:7\:3\:1\in{\cal S}_{10}(1243,2143)$ we obtain 
the essential set ${\cal E}(\pi)=\{(9,1),(8,3),(5,3),(4,4),(4,7),(2,8)\}$ where 
$(4,7)$ is of rank $1$. Replacing this element yields the essential set of 
$\sigma=7\:9\:8\:5\:4\:2\:6\:10\:3\:1\in{\cal S}_{10}(1243,2143)$: 
\begin{center}
\unitlength=0.25cm
\begin{picture}(27,10)
\linethickness{0.3pt}
\multiput(0,0)(0,1){11}{\line(1,0){10}}
\multiput(0,0)(1,0){11}{\line(0,1){10}}
\put(3.5,6.5){\makebox(0,0)[cc]{\sf\tiny0}}
\put(2.5,5.5){\makebox(0,0)[cc]{\sf\tiny0}}
\put(0.5,1.5){\makebox(0,0)[cc]{\sf\tiny0}}
\put(7.5,8.5){\makebox(0,0)[cc]{\sf\tiny1}}
\put(6.5,6.5){\makebox(0,0)[cc]{\sf\tiny1}}
\put(2.5,2.5){\makebox(0,0)[cc]{\sf\tiny1}}
\multiput(17,0)(0,1){11}{\line(1,0){10}}
\multiput(17,0)(1,0){11}{\line(0,1){10}}
\put(20.5,6.5){\makebox(0,0)[cc]{\sf\tiny0}}
\put(19.5,5.5){\makebox(0,0)[cc]{\sf\tiny0}}
\put(17.5,1.5){\makebox(0,0)[cc]{\sf\tiny0}}
\put(24.5,8.5){\makebox(0,0)[cc]{\sf\tiny1}}
\put(22.5,7.5){\makebox(0,0)[cc]{\sf\tiny0}}
\put(19.5,2.5){\makebox(0,0)[cc]{\sf\tiny1}}
\put(4.5,9.5){\circle*{0.5}}
\put(9.5,6.5){\circle*{0.5}}
\put(6.5,2.5){\circle*{0.5}}
\put(23.5,9.5){\circle*{0.5}}
\put(21.5,6.5){\circle*{0.5}}
\put(26.5,2.5){\circle*{0.5}}
\linethickness{0.1pt}
\multiput(4,9)(0,0.2){5}{\line(1,0){6}}
\multiput(4,8)(0,0.2){5}{\line(1,0){1}}\multiput(8,8)(0,0.2){5}{\line(1,0){2}}
\multiput(4,7)(0,0.2){5}{\line(1,0){1}}\multiput(7,7)(0,0.2){5}{\line(1,0){3}}
\multiput(4,6)(0,0.2){5}{\line(1,0){1}}\multiput(7,6)(0,0.2){5}{\line(1,0){3}}
\multiput(3,5)(0,0.2){5}{\line(1,0){7}}
\multiput(1,4)(0,0.2){5}{\line(1,0){9}}
\multiput(1,3)(0,0.2){5}{\line(1,0){1}}\multiput(3,3)(0,0.2){5}{\line(1,0){7}}
\multiput(1,2)(0,0.2){5}{\line(1,0){1}}\multiput(3,2)(0,0.2){5}{\line(1,0){7}}
\multiput(1,1)(0,0.2){5}{\line(1,0){9}}
\multiput(0,0)(0,0.2){5}{\line(1,0){10}}
\multiput(23,9)(0,0.2){5}{\line(1,0){4}}
\multiput(23,8)(0,0.2){5}{\line(1,0){1}}\multiput(25,8)(0,0.2){5}{\line(1,0){2}}
\multiput(23,7)(0,0.2){5}{\line(1,0){4}}
\multiput(21,6)(0,0.2){5}{\line(1,0){6}}
\multiput(20,5)(0,0.2){5}{\line(1,0){7}}
\multiput(18,4)(0,0.2){5}{\line(1,0){9}}
\multiput(18,3)(0,0.2){5}{\line(1,0){1}}\multiput(20,3)(0,0.2){5}{\line(1,0){7}}
\multiput(18,2)(0,0.2){5}{\line(1,0){1}}\multiput(20,2)(0,0.2){5}{\line(1,0){7}}
\multiput(18,1)(0,0.2){5}{\line(1,0){9}}
\multiput(17,0)(0,0.2){5}{\line(1,0){10}}
\put(13.5,5){\makebox(0,0)[cc]{$\longrightarrow$}}
\linethickness{0.7pt}
\multiput(6,6)(0,1){2}{\line(1,0){1}}\multiput(6,6)(1,0){2}{\line(0,1){1}}
\multiput(22,7)(0,1){2}{\line(1,0){1}}\multiput(22,7)(1,0){2}{\line(0,1){1}}
\end{picture}
\vspace*{1ex}

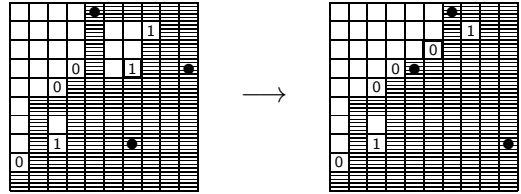
{\footnotesize{\bf Figure 5}\hspace*{0.25cm}\parbox[t]{12cm}{On the left the diagram of $\pi$; 
on the right the diagram of $\sigma$. Only the given dots change their 
position.}}
\end{center}
\erm
\end{bsp}
\vspace*{1ex}

The pattern considered now is closely related to the increasing subsequences. 
In special case of $132$-avoiding permutations the following 
characterization is identical with \cite[Theo. 4.1c]{reifegerste}. 

\begin{satz} \label{213...k-avoiding}
Let $\pi\in{\cal S}_n(1243,2143)$ be a Schr\"oder permutation. 
Then $\pi$ avoids $213\cdots k$ if and only if every element $(i,j)\in{\cal 
E}(\pi)$ satisfies $i+j\ge n+3-k+\rho(i,j)$.
\end{satz}

\begin{bew}
Let $(i,j)\in{\cal E}(\pi)$. The $n$ dots representing $\pi$ are arranged as 
follows (the labels are the numbers of dots contained in the certain regions 
where $r:=\rho(i,j)$):
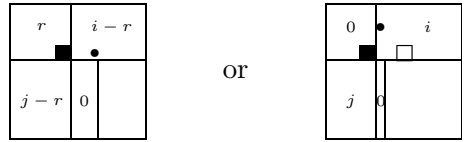
\begin{center}
\unitlength=0.3cm
\begin{picture}(20,6)
\linethickness{0.3pt}
\multiput(0,0)(0,6){2}{\line(1,0){6}}
\multiput(0,0)(6,0){2}{\line(0,1){6}}
\multiput(14,0)(0,6){2}{\line(1,0){6}}
\multiput(14,0)(6,0){2}{\line(0,1){6}}
\put(2,3.5){\rule{0.2cm}{0.2cm}}\put(15.5,3.5){\rule{0.2cm}{0.2cm}}
\put(17,3.4){$\Box$}
\put(0,3.5){\line(1,0){6}}\put(2.7,0){\line(0,1){6}}\put(3.9,0){\line(0,1){3.5}}
\put(14,3.5){\line(1,0){6}}\put(16.2,0){\line(0,1){6}}\put(16.6,0){\line(0,1){3.5}}
\put(10,3){\makebox(0,0)[cc]{or}}
\put(3.75,3.85){\circle*{0.4}}\put(16.42,5){\circle*{0.4}}
\put(1.4,5){\makebox(0,0)[cc]{\tiny$r$}}
\put(4.5,5){\makebox(0,0)[cc]{\tiny$i-r$}}
\put(1.4,1.7){\makebox(0,0)[cc]{\tiny$j-r$}}
\put(3.3,1.7){\makebox(0,0)[cc]{\tiny$0$}}
\put(15.1,5){\makebox(0,0)[cc]{\tiny$0$}}
\put(18.5,5){\makebox(0,0)[cc]{\tiny$i$}}
\put(15.1,1.7){\makebox(0,0)[cc]{\tiny$j$}}
\put(16.435,1.7){\makebox(0,0)[cc]{\tiny$0$}}
\end{picture}
\end{center}
If there is no corner $(i,j')$ such that $j<j'$ (see the left-hand picture) then for all elements $\pi_k$ 
with $k>i$ and $\pi_k>j$ we have $\pi_k>\pi_i$. Clearly, $\pi_{i+1}<\pi_i$. On the other hand, if there exists such a corner $(i,j')$ (see 
the right-hand picture; necessarily, $\rho(i,j)=0$ and $\rho(i,j')=1$) then the dot northwest 
of $(i,j')$ is contained in column $j+1$, otherwise Lemma \ref{diagram 
properties}b fails to hold. (Note that this dot marks an inner corner of the 
$0$-component.) In this case, all elements $\pi_k$ 
with $k>i$ and $\pi_k>j$ satisfy $\pi_k>j+1$. Clearly, $\pi_{i+1}\le j$. In both cases the elements representing by dots in the 
lower right-hand region appear in increasing order since $\pi$ is $2143$-avoiding. If 
$i+j<n+3-k+\rho(i,j)$ then their number $n-(i+j)+\rho(i,j)$ is at least $k-2$.\\
To prove the converse, suppose that every element of the essential set satisfies the 
above condition. Then we have $\pi_i+i>n+3-k$ for all $i\in{\sf D}(\pi)$. 
Hence for each descent $i$ of $\pi$ there exist at most $k-3$ elements $\pi_j$ 
with $j>i$ and $\pi_j>\pi_i$. Since $\pi$ is $2143$-avoiding these elements form an increasing 
sequence. Thus there is no pattern $2134\cdots k$ in $\pi$.
\end{bew}
\vspace*{0.1cm}
    
Now we are in the position to answer the first of a collection of open 
problems given in \cite{egge-mansour}. Theorem 6.5 of this reference implies 
the Wilf-equivalence of $\{1243,2143,12\cdots k\}$ and $\{1243,2143,213\cdots 
k\}$, that is, $|{\cal S}_n(1243,2143,12\cdots k)|=|{\cal S}_n(1243,2143,213\cdots k)|$ 
for all $n$ and $k$. The authors asked for a combinatorial proof of this fact. There is a simple bijection in terms of diagrams 
whose essence was already used to prove the analogue for $132$-avoiding 
permutations.

\begin{kor} \label{bijection omega}
There is a bijection $\omega:{\cal S}_n(1243,2143)\to{\cal S}_n(1243,2143)$ 
such that for all $k\ge 1$ and any permutation $\pi\in{\cal S}_n(1243,2143)$, 
we have that $\pi$ avoids $12\cdots k$ if and only if $\omega(\pi)$ avoids 
$213\cdots k$.
\end{kor}

\begin{bew}
Let $\pi\in{\cal S}_n(1243,2143)$ be a Schr\"oder permutation which avoids 
$12\cdots k$. By Theorem \ref{12...k-avoiding} and \cite[Theo. 4.1b]{reifegerste}, the diagram of $\phi(\pi)\in{\cal S}_n(132)$ contains 
$(n+1-k,n-k,\ldots,1)$. Since all the corners of $(n+1-k,n-k,\ldots,1)$ are in 
the diagonal $i+j=n+2-k$ we have $i+j\ge n+2-k+2\rho(i,j)$ for all 
$(i,j)\in{\cal E}(\pi)$. Hence the diagram corners of rank $1$ satisfy the 
condition of Theorem \ref{213...k-avoiding} anyway. Thus every diagram corresponding 
to a $12\cdots k$-avoiding Schr\"oder permutation is uniquely determined by its 
corners outside the shape $(n+1-k,n-k,\ldots,1)$, that is, by all corners except 
for those satisfying $i+j=n+2-k$. Consequently, the diagram of $\omega(\pi)$ we define to be 
this one whose corners are the corners of $D(\pi)$ which are not contained in 
$(n+1-k,n-k,\ldots,1)$. Ranks are kept up.\\
Conversely, given any Schr\"oder permutation $\sigma\in{\cal S}_n(1243,2143)$ 
whose all diagram corners satisfy $i+j\ge n+3-k+\rho(i,j)$ we construct the 
permutation $\omega^{-1}(\sigma)$ as follows: let $E$ be the corner set 
of the diagram obtained as union of 
$D(\phi(\sigma))$ and $(n+1-k,n-k,\ldots,1)$. (Note that this is a Young 
diagram since $D(\phi(\pi))$ is such a one.) Then we form the essential set of 
$\omega^{-1}(\sigma)$ from the pairs $(i,j)\in E$ 
for which $(i+1,j+1)\notin{\cal E}_1(\sigma)$, and all elements of ${\cal 
E}_1(\sigma)$. The first ones are defined to be of rank $0$, the rank of the latter 
should be $1$. Obviously, the resulting set is an essential set of a $12\cdots 
k$-avoiding Schr\"oder permutation. 
\end{bew}

\begin{bsp}
\brm
The maximum length of an increasing subsequence in the Schr\"oder permutation $\pi=4\:6\:3\:1\:5\:7\:2\in{\cal S}_7(1243,2143)$ equals $3$. Taking $k=4$, we 
obtain $\omega(\pi)=1\:6\:3\:4\:5\:7\:2$:
\begin{center}
\unitlength=0.3cm
\begin{picture}(20,7)
\linethickness{0.3pt}
\multiput(0,0)(0,1){8}{\line(1,0){7}}
\multiput(0,0)(1,0){8}{\line(0,1){7}}
\multiput(13,0)(0,1){8}{\line(1,0){7}}
\multiput(13,0)(1,0){8}{\line(0,1){7}}
\put(4.5,5.5){\makebox(0,0)[cc]{\sf\tiny1}}
\put(2.5,5.5){\makebox(0,0)[cc]{\sf\tiny0}}
\put(1.5,4.5){\makebox(0,0)[cc]{\sf\tiny0}}
\put(1.5,1.5){\makebox(0,0)[cc]{\sf\tiny1}}
\put(17.5,5.5){\makebox(0,0)[cc]{\sf\tiny1}}
\put(14.5,1.5){\makebox(0,0)[cc]{\sf\tiny1}}
\bezier{50}(1,4)(1.5,4.5)(2,5)
\bezier{50}(1,5)(1.5,4.5)(2,4)
\bezier{50}(2,5)(2.5,5.5)(3,6)
\bezier{50}(2,6)(2.5,5.5)(3,5)
\linethickness{0.1pt}
\multiput(3,6)(0,0.2){5}{\line(1,0){4}}
\multiput(3,5)(0,0.2){5}{\line(1,0){1}}\multiput(5,5)(0,0.2){5}{\line(1,0){2}}
\multiput(2,4)(0,0.2){5}{\line(1,0){5}}
\multiput(0,3)(0,0.2){5}{\line(1,0){7}}
\multiput(0,2)(0,0.2){5}{\line(1,0){1}}\multiput(2,2)(0,0.2){5}{\line(1,0){5}}
\multiput(0,1)(0,0.2){5}{\line(1,0){1}}\multiput(2,1)(0,0.2){5}{\line(1,0){5}}
\multiput(0,0)(0,0.2){5}{\line(1,0){7}}
\multiput(13,6)(0,0.2){5}{\line(1,0){7}}
\multiput(13,5)(0,0.2){5}{\line(1,0){1}}\multiput(18,5)(0,0.2){5}{\line(1,0){2}}
\multiput(13,4)(0,0.2){5}{\line(1,0){1}}\multiput(15,4)(0,0.2){5}{\line(1,0){5}}
\multiput(13,3)(0,0.2){5}{\line(1,0){1}}\multiput(15,3)(0,0.2){5}{\line(1,0){5}}
\multiput(13,2)(0,0.2){5}{\line(1,0){1}}\multiput(15,2)(0,0.2){5}{\line(1,0){5}}
\multiput(13,1)(0,0.2){5}{\line(1,0){1}}\multiput(15,1)(0,0.2){5}{\line(1,0){5}}
\multiput(13,0)(0,0.2){5}{\line(1,0){7}}
\linethickness{0.7pt}
\multiput(4,5)(0,1){2}{\line(1,0){1}}\multiput(4,5)(1,0){2}{\line(0,1){1}}
\multiput(2,5)(0,1){2}{\line(1,0){1}}\multiput(2,5)(1,0){2}{\line(0,1){1}}
\multiput(1,4)(0,1){2}{\line(1,0){1}}\multiput(1,4)(1,0){2}{\line(0,1){1}}
\multiput(1,1)(0,1){2}{\line(1,0){1}}\multiput(1,1)(1,0){2}{\line(0,1){1}}
\multiput(17,5)(0,1){2}{\line(1,0){1}}\multiput(17,5)(1,0){2}{\line(0,1){1}}
\multiput(14,1)(0,1){2}{\line(1,0){1}}\multiput(14,1)(1,0){2}{\line(0,1){1}}
\put(10,3.5){\makebox(0,0)[cc]{$\longleftrightarrow$}}
\end{picture}
\vspace*{1ex}

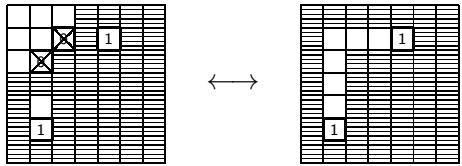
{\footnotesize{\bf Figure 6}\hspace*{0.25cm}Construction of $\omega(\pi)$: the 
corners crossed out satisfy $i+j=n+2-k$.}
\end{center}
\erm
\end{bsp}

\begin{beme}
\brm
\begin{enum}
\item[]
\item Since ${\cal E}_1(\pi)={\cal E}_1(\omega(\pi))$ for all $\pi\in{\cal 
S}(1243,2143,12\cdots k)$ the map $\omega$ takes any $132$-avoi\-ding permutation to 
a permutation which avoids $132$ as well. Indeed, the restriction of $\omega$ on 
${\cal S}(132,12\cdots k)$ is precisely the bijection given in 
\cite[Cor. 4.4]{reifegerste} that proves the Wilf-equivalence of $\{132,12\cdots k\}$ and 
$\{132,213\cdots k\}$.
\item It is clear from the construction that a $12\cdots k$-avoiding Schr\"oder permutation 
also avoids $213\cdots k$ if and only if it is a fixed point of 
$\omega$. The essential set of such a permutation can be constructed as follows: 
consider the corner set of a Young diagram which contains 
$(n+1-k,n-k,\ldots,1)$, and fits in $(n-1,n-2,\ldots,1)$. Now replace at least 
all elements $(i,j)$ by $(i+1,j+1)$ for which $i+j=n+2-k$. Some further corners can be 
replaced if these satisfy $i+j<n$. The rank of all new corners is set as $1$; 
let the others be of rank $0$.\\
We will discuss the enumerative consequence only for $k=3$.
\end{enum}
\erm
\end{beme}

\begin{kor}
$|{\cal S}_n(1243,2143,123,213)|=2^{n-1}$ for all $n\ge 1$.
\end{kor}

\begin{bew}
Taking up again the idea of the previous remark, the diagram of a Schr\"oder 
permutation $\pi\in{\cal S}_n(1243,2143)$ which avoids both $123$ as $213$ 
arises from a Young diagram that contains $(n-2,n-3,\ldots,1)$, and fits in $(n-1,n-2,\ldots,1)$. 
Clearly, each such Young diagram is uniquely determined by its corners in the 
diagonal $i+j=n$. In particular, there are $2^{n-1}$ diagrams of this kind. 
(This implies $|{\cal S}_n(132,123)|=2^{n-1}$; see \cite[Prop. 
7]{simion-schmidt} for another proof.)\\ 
From the corner set of each Young diagram the essential set of only one permutation $\pi\in{\cal S}_n(1243,2143,123,213)$ 
can be generated because all the corners $(i,n-1-i)$ must be replaced, but all 
the corners $(i,n-i)$ must not be replaced.
\end{bew}

\begin{bem}
\brm
To obtain a $\{132,123,213\}$-avoiding permutation in the way described in the 
proof, the Young diagram must not have any corner in the diagonal $i+j=n-1$. We 
can identify such a diagram by a binary sequence of length $n-1$ whose $i$th 
element is defined as $1$ (or $0$) if $(i,n-i)$ is a corner (or not). The 
condition that there is no diagram corner outside the diagonal $i+j=n$ means 
that the corresponding sequence contains no consecutive zeros. The number of 
such sequences is known to be equal the $(n+1)$st Fibonacci number $F_{n+1}$. (The {\it Fibonacci numbers} are defined 
by $F_1=F_2=1$ and $F_n=F_{n-1}+F_{n-2}$ for $n\ge3$.) The result $|{\cal 
S}_n(132,123,213)|=F_{n+1}$ already appears in \cite[Prop. 15]{simion-schmidt}.
\erm
\end{bem}

The next result deals with the occurrence of a decreasing subsequence of 
length $k$ in Schr\"oder permutations. The analogue for $132$-avoiding 
permutations is simple: a permutation $\pi\in{\cal S}_n(132)$ avoids 
$k(k-1)\cdots1$ if and only if $|{\cal E}(\pi)|\le k-2$, see \cite[Theo. 
4.1a]{reifegerste}. Now the condition is some more difficult.\\
To state it, we first set some notation. For a permutation $\pi\in{\cal S}_n$ 
we denote by $r(\pi)$ and $c(\pi)$ the number of rows and 
columns, respectively, that contain a diagram corner. As mentioned above, we 
have $r(\pi)=\des(\pi)$, and $c(\pi)=\des(\pi^{-1})$. (Note that the transpose 
of $D(\pi)$ is just the diagram of $\pi^{-1}$.) It follows from Proposition 
\ref{fulton characterization} that any diagram row (column) contains at most two corners 
(necessarily of different rank) if $\pi\in{\cal S}_n(1243,2143)$. Let $r_2(\pi)$ 
and $c_2(\pi)$ be the number of diagram rows and diagram columns, respectively, containing two corners.

\begin{satz} \label{k...1-avoiding}
Let $\pi\in{\cal S}_n(1243,2143)$ be a Schr\"oder permutation. 
Then $\pi$ avoids $k(k-1)\cdots1$ if and only if one of the following 
conditions holds:
\begin{enum2}
\item $r(\pi)\le k-2$ or $c(\pi)\le k-2$; 
\item $r(\pi)=k-1,\;r_2(\pi)=c_2(\pi)=1$, and there is no element 
$(i,j)\in{\cal E}(\pi)$ that such both row $i$ and column $j$ contain another 
corner.
\end{enum2}
\end{satz}
\vspace*{-0.35cm}

\begin{bew}
Suppose that $\pi$ contains a decreasing subsequence of length $k$. Obviously, 
its inverse contains such a sequence as well. Consequently, both $\pi$ and $\pi^{-1}$ must have at least $k-1$ descents, that is, $r(\pi)\ge k-1$ and 
$c(\pi)\ge k-1$. Now let $r(\pi)=k-1$ and $r_2(\pi)=c_2(\pi)=1$. (Then we also 
have $c(\pi)=k-1$.)\\
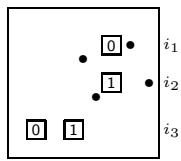
\begin{minipage}[t]{16cm}
\raisebox{-1.9cm}
{\unitlength=0.25cm
\begin{picture}(8,8.5)
\multiput(0,0)(0,8){2}{\line(1,0){8}}
\multiput(0,0)(8,0){2}{\line(0,1){8}}
\put(6.5,6){\circle*{0.4}}
\put(7.5,4){\circle*{0.4}}
\put(4,5.25){\circle*{0.4}}
\put(4.7,3.25){\circle*{0.4}}
\multiput(5,5.5)(0,1){2}{\line(1,0){1}}\multiput(5,5.5)(1,0){2}{\line(0,1){1}}
\multiput(5,3.5)(0,1){2}{\line(1,0){1}}\multiput(5,3.5)(1,0){2}{\line(0,1){1}}
\multiput(1,1)(0,1){2}{\line(1,0){1}}\multiput(1,1)(1,0){2}{\line(0,1){1}}
\multiput(3,1)(0,1){2}{\line(1,0){1}}\multiput(3,1)(1,0){2}{\line(0,1){1}}
\put(5.5,6){\makebox(0,0)[cc]{\sf\tiny0}}
\put(5.5,4){\makebox(0,0)[cc]{\sf\tiny1}}
\put(1.5,1.5){\makebox(0,0)[cc]{\sf\tiny0}}
\put(3.5,1.5){\makebox(0,0)[cc]{\sf\tiny1}}
\put(8.25,6){\makebox(0,0)[lc]{\tiny$i_1$}}
\put(8.25,4){\makebox(0,0)[lc]{\tiny$i_2$}}
\put(8.25,1.5){\makebox(0,0)[lc]{\tiny$i_3$}}
\end{picture}}
\hfill
\parbox[t]{13cm}
{If the essential corners are arrange as in the picture opposite (where 
$i_2\not=i_3$) we have $\pi_{i_1}>\pi_{i_1+1},\;\pi_{i_2}>\pi_{i_2+1}$ (corners 
correspond to descents) but $\pi_{i_1}<\pi_{i_2}$, and 
$\pi_{i_1+1}<\pi_{i_2+1}$. The last relation follows from the fact that 
$(i_2,\pi_{i_2}+1)$ is a diagram square by the construction. If its rank would 
be equal to $0$ then} 
\end{minipage}
\vspace*{-0.05cm} 

$r_2>1$. Since $\des(\pi)=k-1$ there is no decreasing 
subsequence of length $k$ in $\pi$ which contradicts the assumption. 
In the second case (corners $(i_1,j_2),\;(i_2,j_1)$ of rank $0$, corners 
$(i_1,j_3),\;(i_3,j_1)$ of rank $1$ where $i_1<i_2<i_3$ and $j_1<j_2<j_3$) use 
the same arguments.\\
On the other hand, if condition (i) holds then $\des(\pi)\le k-2$ or 
$\des(\pi^{-1})\le k-2$ and hence $\pi\in{\cal S}_n(k\cdots1)$. If (ii) is 
satisfied then (as shown in the first part of the proof) $\pi$ cannot contain any decreasing subsequence of length $k$.  
\end{bew}
\vspace*{1ex}

Here we will enumerate the permutations described in Theorem 
\ref{k...1-avoiding} only for $k=3$. To satisfy condition (ii) is impossible 
in this case. Thus a Schr\"oder permutation is $321$-avoiding if and only if 
all its diagram corners are either in the same row or in the same column. This 
characterization was already given in \cite[Prop. 5.4]{eriksson-linusson} for 
$321$-avoiding vexillary permutations. (Note that the essential set of a 
vexillary permutation can contain elements of rank greater than $1$; for 
example, $1243$ is such a permutation.)\\ 
Egge and Mansour have shown (derived from the generating function in 
\cite[Prop. 7.4]{egge-mansour}) that
\bdpm
|{\cal S}_n(1243,2143,321)|={n-1\choose 0}+{n-1\choose 1}+2{n-1\choose 2}+2{n-1\choose 
3}\quad\mbox{for all }n\ge1.
\edpm
Their fourth problem asked for a combinatorial proof. Here it is.

\begin{kor}
$|{\cal S}_n(1243,2143,321)|=n+2{n\choose3}$ for all $n\ge 1$.
\end{kor}

\begin{bew}
Let $\pi$ be a Schr\"oder permutation avoiding $321$. We distinguish the three cases mentioned at the begin of the section.\\ 
If $\pi\in{\cal S}_n(132)$, different from the identity, then its diagram is 
a rectangle whose lower right-hand corner $(i,j)$ satisfies $i+j\le n$. There 
are 
\bdpm
\sum_{i=1}^{n-1} (n-i)={n\choose2}
\edpm
such diagrams. (The enumeration of $\{132,321\}$-avoiding permutations was first done in \cite[Prop. 11]{simion-schmidt}.) If there exists no element of 
rank $0$ in ${\cal E}(\pi)$, the permutation $\pi_2-1\cdots\pi_n-1$ belongs to 
${\cal S}_{n-1}(132,321)$.\\
It remains to consider the case that $D(\pi)$ has corners of rank $0$ and $1$. Since these squares are in the same row or column there is exactly one 
corner of each rank. Without loss of generality, we may assume that the both elements 
of ${\cal E}(\pi)$ are in the same row. (For the result in terms 
of columns consider the transpose that corresponds to the inverse of $\pi$.) 
Let $(i,j)\in{\cal E}_0(\pi)$ and $(i,j')\in{\cal E}_1(\pi)$. From Proposition \ref{fulton characterization} 
the conditions $1<i$, $j+1<j'$, and $i+j'\le n+1$ results.
Given the pair $(i,j)$, the integer $j'$ can be 
chosen in $n-i-j$ ways where $i\in\{2,\ldots,n-2\}$, and $j\in\{1,\ldots,n-1-i\}$. Clearly,
\bdpm
\sum_{i=2}^{n-2}\sum_{j=1}^{n-1-i} (n-i-j)=\frac{1}{2}\sum_{i=1}^{n-3} 
i(i+1)={n-1\choose3}.
\edpm
Summarized, we obtain $|{\cal 
S}_n(1243,2143,321)|=1+{n\choose2}+{n-1\choose2}+2{n-1\choose3}$. (Note that 
the term $1$ stands for the identity; the factor $2$ regards rows and columns 
in the third case.)
\end{bew}
\vspace*{1ex}

The last pattern we will discuss is a special case of an important class as 
well. In \cite[Theo. 4.5]{reifegerste}, we characterized $132$-avoiding 
permutations which avoid the additional pattern $s(s+1)\cdots k12\cdots(s-1)$ 
where $s\in\{2,\ldots,k\}$, and $k\ge3$. The condition given there was of 
technical nature but for $s=2$ and $k=3$ it is equivalent to the following 
simple one: a permutation $\pi\in{\cal S}_n(132)$ avoids $231$ if and only if 
all its diagram rows are of distinct length, that means, all diagram rows 
contain a corner. 
Analogously to that, $231$-avoiding Schr\"oder permutations can be described.

\begin{prop} \label{231-avoiding}
A Schr\"oder permutation $\pi\in{\cal S}_n(1243,2143)$ avoids $231$ if and 
only if 
\begin{enum2}
\item every diagram row contains a (namely, exactly one) element of the essential set, 
\item and every diagram column contains at most an element of the essential set. 
\end{enum2}
\end{prop}
\vspace*{-0.35cm}

\begin{bew}
Suppose that $\pi\in{\cal S}_n(1243,2143)$ contains a $231$-pattern: let $i_1<i_2<i_3$ 
such that $\pi_{i_3}<\pi_{i_1}<\pi_{i_2}$. 
Since diagram corners and permutation descents corresponds to each other, there 
is an integer $i$ with $i_1\le i<i_2$ such that row $i$ contains no corner. We 
assume that this row does not belong to the diagram, otherwise condition (i) fails to 
hold. Since $(i_2,\pi_{i_3})\in D(\pi)$ we have $\pi_i=1$, and all diagram 
squares appearing below the $i$th row are of rank $1$. By the construction, 
$(i+1,2)$ is the upper left-hand corner of the component which contains 
$(i_2,\pi_{i_3})$ (and all the other diagram squares of rank $1$). By Lemma 
\ref{diagram properties}b there is no corner in the strict northwest of another 
one. Therefore, and since $\pi_{i_1}<\pi_{i_2}$, the diagram corners containing 
in row $i-1$ and $i+1$, respectively, have to be in the same column.\\
For the other direction, we suppose that there are a diagram square $(i_1,j)$ 
for which $(i_1,j+1)\notin D(\pi)$, and a square $(i_2,j)\in{\cal E}(\pi)$ with 
$i_1<i_2$. (Then one of the conditions (i) and (ii) is not satisfied.) 
It is easy to see that $\pi_{i_1}\pi_{i_2}\pi_{i_3}$ where $\pi_{i_3}=j$ is a 
subsequence of type $231$.
\end{bew}

From the first part of the proof, it is clear how the diagram of a 
$231$-avoiding Schr\"oder permutation has to look. 

\begin{kor} \label{231-avoiding corollary}
The diagram of a Schr\"oder permutation satisfies the conditions of {\rm 
Proposition \ref{231-avoiding}} if and only if it 
is of the following shape:
\begin{center}
\unitlength=0.25cm
\begin{picture}(10,10)
\linethickness{0.2pt}
\multiput(0,0)(0,10){2}{\line(1,0){10}}
\multiput(0,0)(10,0){2}{\line(0,1){10}}
\put(0.5,4.5){\circle*{0.5}}
\put(3,8){\makebox(0,0)[cc]{\tiny$E_0$}}
\put(2,3.25){\makebox(0,0)[cc]{\tiny$E_1$}}
\linethickness{0.7pt}
\put(0,5){\line(1,0){5}}\put(5,5){\line(0,1){1}}
\put(5,6){\line(1,0){1}}\put(6,6){\line(0,1){1}}
\put(9,9){\line(0,1){1}}\put(8,9){\line(1,0){1}}
\bezier{3}(6.5,7.5)(7.25,8.25)(8,9)
\put(1,4){\line(1,0){3}}\put(4,3){\line(0,1){1}}
\put(3,3){\line(1,0){1}}\put(2,1){\line(0,1){1}}
\put(1,1){\line(1,0){1}}\put(1,1){\line(0,1){3}}
\put(0,5){\line(0,1){5}}\put(0,10){\line(1,0){9}}
\bezier{3}(6.5,7.5)(7.25,8.25)(8,9)
\put(2.5,2.5){\makebox(0,0)[cc]{$\cdot$}}
\end{picture}
\end{center}
where the diagram components $E_0$ and $E_1$ are Young diagrams whose row 
lengths are each distinct. (It may be that $E_0$ and/or $E_1$ are empty.) 
\end{kor}

As an immediate consequence we can characterize $231$-avoiding Schr\"oder 
permutations from their descent set. 

\begin{kor} 
Let $\pi\in{\cal S}_n$ be a permutation, and ${\sf D}(\pi)$ the set of its 
descents. Then $\pi$ avoids the patterns $1243$, $2143$, and $231$ if and 
only if one of the following conditions is satisfied:
\begin{enum2}
\item $s>d$ and ${\sf D}(\pi)=\{1,2,\ldots,d\}$;
\item $s\le d$ and ${\sf D}(\pi)=\{1,2,\ldots,s-1,s+1,s+2,\ldots,d+1\}$, and 
$\pi_{s-1}>\pi_{s+1}$ if $1<s<n$
\end{enum2}
\vspace*{-0.35cm}

where $\pi_s=1$ and $d=\des(\pi)$. 
\end{kor}

\begin{bem}
\brm
We can decide by the conditions whether a $231$-avoiding Schr\"oder permutation 
avoids the pattern $132$ in addition or not. The first condition describes all 
permutations in ${\cal S}_n(132,231)$.  
\erm
\end{bem}

Corollary \ref{231-avoiding corollary} yields the answer to the third question of the Egge-Mansour list.

\begin{kor}
For $n\ge2$ we have $|{\cal S}_n(1243,2143,231)|=(n+2)2^{n-3}$.
\end{kor}

\begin{bew}
Let $\pi\in{\cal S}_n(1243,2143,231)$. Consider the partition $\lambda(\pi)$ 
whose parts are just the lengths of the diagram rows where the length 
of the first row containing squares of rank 1 is listed twice. (By this 
information $D(\pi)$ and hence $\pi$ is completely described.) Adding some zeros (if 
necessary), we may assume that 
$\lambda=(\lambda_1,\ldots,\lambda_{n-1})$ is of length $n-1$. 
By the previous discussion, 
we have 
\bdpm
n>\lambda_1>\lambda_2>\ldots>\lambda_{i-1}\ge\lambda_i>\lambda_{i+1}>\ldots>
\lambda_l>\lambda_{l+1}=\lambda_{l+2}=\ldots=\lambda_{n-1}=0
\edpm 
with $l\in\{0,\ldots,n-1\}$, and $\lambda\subseteq(n-1,n-1,n-2,\ldots,2)$ 
(where $\subseteq$ means the containment of the corresponding diagrams). 
Furthermore we have $\lambda_j\not=1$ for all $j$ if any positive part exists 
twice.\\
To set $\lambda$ such that all the positive parts of $\lambda$ are distinct, 
there are $\sum_{l=0}^{n-1} {n-1\choose l}=2^{n-1}$ ways. (This case corresponds to 132-avoiding permutations; hence $|{\cal 
S}_n(132,231)|=2^{n-1}$, see also \cite[Prop. 9]{simion-schmidt}.) If there 
exists a positive part twice then the number of partitions be considered equals 
\bdpm
\sum_{l=1}^{n-2} {n-2\choose l}{l\choose 1}=(n-2)\sum_{l=0}^{n-3} {n-3\choose 
l}=(n-2)2^{n-3}.
\edpm
Note that $1$ must not a part of $\lambda$ now. Consequently, there are 
$2^{n-1}+(n-2)2^{n-3}=(n+2)2^{n-3}$ Schr\"oder permutations in ${\cal S}_n$ 
which avoid $231$.
\end{bew}
\vspace*{0.75cm}


\setcounter{section}{4}\setcounter{satz}{0}

\centerline{\large{\bf 4}\hspace*{0.25cm}
{\sc A correspondence to lattice paths}}
\vspace*{0.5cm}

It is well known that the $n$th Schr\"oder number $r_n$ counts the number of all lattice 
paths from the origin to $(n,n)$, with steps $[1,0]$ (called {\sf E}ast steps), $[0,1]$ (called 
{\sf N}orth steps), and $[1,1]$ (called {\sf D}iagonal steps), that never pass below the 
line $y=x$. Such paths we call {\it Schr\"oder paths}. (See \cite[Exc. 6.39]{stanley} for further combinatorial interpretations of 
the Schr\"oder numbers.)\\[2ex]
Egge and Mansour have given a bijection $\Psi_{EM}$ between these paths and 
Schr\"oder permutations in ${\cal S}_{n+1}$, see \cite[Sect. 4]{egge-mansour}. 
Its essential property is: the number of subsequences $12\cdots k$ occuring in 
$\pi\in{\cal S}_{n+1}(1243,2143)$ can read off (more or less) directly 
from the path $\Psi_{EM}(\pi)$.\\ 
This bijection can be unterstood as analogue of Krattenthaler's correspondence 
$\Psi_K$ between $132$-avoiding permutations in ${\cal S}_n$ and lattice paths 
from $(0,0)$ to $(n,n)$ without diagonal steps, never passing below the line $y=x$. The map $\Psi_K$ 
encodes the number of increasing subsequences of prescribed length in the same 
way, see \cite[(3.2)]{krattenthaler}.\\[2ex] 
We pointed out in \cite{reifegerste} that the path $\Psi_K(\pi)$ and the 
diagram of a $132$-avoiding permutation $\pi$ are closely related to each 
other. Considering the diagram of $\pi\in{\cal S}_n(132)$ as being contained in 
an $n\times n$-rectangle, $\Psi_K(\pi)$ is the lattice path which goes from the 
upper right-hand to the lower left-hand corners of the rectangle, and travels 
along the diagram boundary.\\
For $\pi=6\:4\:5\:3\:2\:7\:1\in{\cal S}_7(132)$, for example, $\Psi_K(\pi)$ equals the lattice 
path printed in bold: 
\begin{center}
\unitlength=0.3cm
\begin{picture}(7,7)
\linethickness{0.3pt}
\multiput(0,0)(0,1){8}{\line(1,0){7}}
\multiput(0,0)(1,0){8}{\line(0,1){7}}
\linethickness{0.1pt}
\multiput(5,6)(0,0.2){5}{\line(1,0){2}}
\multiput(3,5)(0,0.2){5}{\line(1,0){4}}
\multiput(3,4)(0,0.2){5}{\line(1,0){4}}
\multiput(2,3)(0,0.2){5}{\line(1,0){5}}
\multiput(1,2)(0,0.2){5}{\line(1,0){6}}
\multiput(1,1)(0,0.2){5}{\line(1,0){6}}
\multiput(0,0)(0,0.2){5}{\line(1,0){7}}
\linethickness{1pt}
\put(5,7){\line(1,0){2}}\put(5,6){\line(0,1){1}}
\put(3,6){\line(1,0){2}}\put(3,4){\line(0,1){2}}
\put(2,4){\line(1,0){1}}\put(2,3){\line(0,1){1}}
\put(1,3){\line(1,0){1}}\put(1,1){\line(0,1){2}}
\put(0,1){\line(1,0){1}}\put(0,0){\line(0,1){1}}
\end{picture}
\vspace*{1ex}

{\footnotesize{\bf Figure 7}\hspace*{0.25cm}Lattice path $\Psi_K(6453271)$.}
\end{center}
We can construct the path $\Psi_{EM}$ just as simple from the permutation 
diagram.\\ 
Note again that each (Schr\"oder) permutation is uniquely determined by its 
ranked essential set. Consequently, position and rank of the diagram 
corners are all that we have to transfer to a path corresponding to the permutation. 
Since a permutation $\pi\in{\cal S}_{n+1}(1243,2143)$ should correspond to a 
Schr\"oder path from $(0,0)$ to $(n,n)$, we cannot use $D(\pi)$ itself but the 
diagram of $\phi(\pi)$ is fitting. Recall that the diagram of $\phi(\pi)$ is 
obtained from that one of $\pi$ by "moving" each diagram square of rank $1$ 
northwestwards. Hence it is a Young diagram which is contained in 
$(n-1,n-2,\ldots,1)$. Labeling all corners with their original rank we have all 
information needed to recover $\pi$.\\ 
Now the Schr\"oder path corresponding to $\pi$ is constructed as follows: let 
$D(\phi(\pi))$ be embeded in an $(n-1)\times(n-1)$-rectangle. Analogously to the 
construction of $\Psi_K$, the lattice path is defined to go from the upper right-hand to the lower left-hand corners 
of the rectangle, along the diagram boundary, where every step sequence {\sf NE} 
representing a corner labeled with $0$ is replaced by a step {\sf D}. Last we 
convert the path into the form used in \cite{egge-mansour}. To this end, 
the rectangle is reflected such that the origin is placed at the bottom left instead at the top right. 

\begin{bsp}
\brm
Let $\pi=4\:7\:5\:2\:6\:3\:1\in{\cal S}_7(1243,2143)$.
Using the diagram of $\phi(\pi)$ we can immediately determine the Schr\"oder 
path corresponding to $\pi$ (printed in bold again):
\begin{center}
\unitlength=0.3cm
\begin{picture}(31,7)
\linethickness{0.3pt}
\multiput(0,0)(0,1){8}{\line(1,0){7}}
\multiput(0,0)(1,0){8}{\line(0,1){7}}
\put(5.5,5.5){\makebox(0,0)[cc]{\sf\tiny1}}
\put(2.5,4.5){\makebox(0,0)[cc]{\sf\tiny0}}
\put(0.5,1.5){\makebox(0,0)[cc]{\sf\tiny0}}
\put(2.5,2.5){\makebox(0,0)[cc]{\sf\tiny1}}
\multiput(13,0.5)(0,1){7}{\line(1,0){6}}
\multiput(13,0.5)(1,0){7}{\line(0,1){6}}
\put(17.5,6){\makebox(0,0)[cc]{\sf\tiny1}}
\put(14.5,3){\makebox(0,0)[cc]{\sf\tiny1}}
\put(13,0){\makebox(0,0)[cc]{\tiny$(6,6)$}}
\put(19,7){\makebox(0,0)[cc]{\tiny$(0,0)$}}
\multiput(25,0.5)(0,1){7}{\line(1,0){6}}
\multiput(25,0.5)(1,0){7}{\line(0,1){6}}
\put(25,0){\makebox(0,0)[cc]{\tiny$(0,0)$}}
\put(31,7){\makebox(0,0)[cc]{\tiny$(6,6)$}}
\linethickness{0.1pt}
\multiput(3,6)(0,0.2){5}{\line(1,0){4}}
\multiput(3,5)(0,0.2){5}{\line(1,0){1}}\multiput(6,5)(0,0.2){5}{\line(1,0){1}}
\multiput(3,4)(0,0.2){5}{\line(1,0){4}}
\multiput(1,3)(0,0.2){5}{\line(1,0){6}}
\multiput(1,2)(0,0.2){5}{\line(1,0){1}}\multiput(3,2)(0,0.2){5}{\line(1,0){4}}
\multiput(1,1)(0,0.2){5}{\line(1,0){6}}
\multiput(0,0)(0,0.2){5}{\line(1,0){7}}
\multiput(18,5.5)(0,0.2){5}{\line(1,0){1}}
\multiput(16,4.5)(0,0.2){5}{\line(1,0){3}}
\multiput(16,3.5)(0,0.2){5}{\line(1,0){3}}
\multiput(15,2.5)(0,0.2){5}{\line(1,0){4}}
\multiput(14,1.5)(0,0.2){5}{\line(1,0){5}}
\multiput(14,0.5)(0,0.2){5}{\line(1,0){5}}
\linethickness{1pt}
\put(18,6.5){\line(1,0){1}}\put(18,5.5){\line(0,1){1}}
\put(16,5.5){\line(1,0){2}}\put(16,4.5){\line(0,1){1}}
\bezier{20}(15,3.6)(15.45,4)(15.9,4.4)\put(15,2.5){\line(0,1){1}}
\put(14,2.5){\line(1,0){1}}\put(14,1.5){\line(0,1){1}}
\bezier{20}(13,0.6)(13.45,1)(13.9,1.4)
\put(10,3.5){\makebox(0,0)[cc]{$\longrightarrow$}}
\put(22,3.5){\makebox(0,0)[cc]{$\longrightarrow$}}
\put(25,0.5){\line(0,1){1}}\put(25,1.5){\line(1,0){1}}
\put(26,1.5){\line(0,1){2}}\put(26,3.5){\line(1,0){1}}
\bezier{30}(27,3.5)(27.5,4)(28,4.5)\put(28,4.5){\line(1,0){1}}
\put(29,4.5){\line(0,1){1}}\put(29,5.5){\line(1,0){1}}
\bezier{30}(30,5.5)(30.5,6)(31,6.5)
\end{picture}
\vspace*{1ex}

{\footnotesize{\bf Figure 8}\hspace*{0.25cm}\parbox[t]{12cm}{
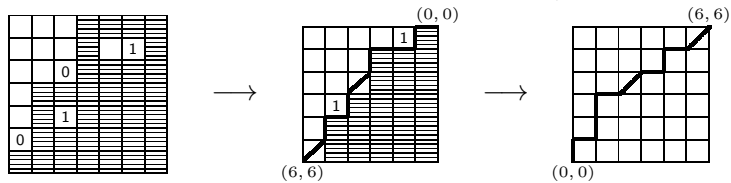
Construction of the path $\Psi_{EM}(\pi)$: On the left the diagram of $\pi$; in 
the centre the diagram of $\phi(\pi)$ with plotted path; on the right the 
converted path.}}
\end{center}
\erm
\end{bsp}

By Lemma \ref{diagram properties}a, each element $(i,j)\in{\cal E}_1(\pi)$ satisfies $i+j\le 
n+1$. Thus, for every corner $(i',j')$ of $D(\phi(\pi))$ labeled with $1$ we have $i'+j'\le n-1$. Therefore, and since $D(\phi(\pi))$ is contained in 
$(n-1,n-2,\ldots,1)$, this construction indeed yields a Schr\"oder path.\\
It is not difficult to see that the path obtained in this way is just $\Psi_{EM}(\pi)$ for $\pi\in{\cal 
S}_n(1243,2143)$ but the construction via diagram requires less effort.\\[2ex]
In \cite{egge-mansour}, the path statistic $\tau_k$ corresponding to the number of 
subsequences of type $12\cdots k$ in $\pi\in{\cal S}_n(1243,2143)$ is defined for $k\ge2$ 
by 
\bdpm
\sum_{s\in\{\sf E,D\}}{h(s)\choose k-1}
\edpm
where $h(s)$ denotes the height of the 
starting point of step $s$. (The {\it height} of a point $(x,y)$ in the 
plane we define to be the difference $y-x$.)
\begin{bsp}
\brm
Consider the Schr\"oder path {\sf NENNEDENED} appearing in the previous example. For each 
east and diagonal step the height is given in the picture.
\begin{center}
\unitlength=0.3cm
\begin{picture}(6,6)
\linethickness{0.1pt}
\multiput(0,0)(0,1){7}{\line(1,0){6}}
\multiput(0,0)(1,0){7}{\line(0,1){6}}
\linethickness{0.7pt}
\put(0,0){\line(0,1){1}}\put(0,1){\line(1,0){1}}
\put(1,1){\line(0,1){2}}\put(1,3){\line(1,0){1}}
\bezier{30}(2,3)(2.5,3.5)(3,4)\put(3,4){\line(1,0){1}}
\put(4,4){\line(0,1){1}}\put(4,5){\line(1,0){1}}
\bezier{30}(5,5)(5.5,5.5)(6,6)
\put(0,1){\circle*{0.2}}
\put(1,3){\circle*{0.2}}
\put(2,3){\circle*{0.2}}
\put(3,4){\circle*{0.2}}
\put(4,5){\circle*{0.2}}
\put(5,5){\circle*{0.2}}
\put(-0.2,1.2){\makebox(0,0)[cb]{\sf\tiny 1}}
\put(0.8,3.2){\makebox(0,0)[cb]{\sf\tiny 2}}
\put(1.8,3.2){\makebox(0,0)[cb]{\sf\tiny 1}}
\put(2.8,4.2){\makebox(0,0)[cb]{\sf\tiny 1}}
\put(3.8,5.2){\makebox(0,0)[cb]{\sf\tiny 1}}
\put(4.8,5.2){\makebox(0,0)[cb]{\sf\tiny 0}}
\end{picture}
\vspace*{1ex}

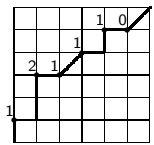
{\footnotesize{\bf Figure 9}\hspace*{0.25cm}Schr\"oder path 
with step heights.}
\end{center}
Thus there occur six subsequences of type $12$ (noninversions), one of type 
$123$, and no one of type $12\cdots k$ for $k\ge4$ in the corresponding permutation $\pi=4\:7\:5\:2\:6\:3\:1\in{\cal S}_7(1243,2143)$ 
\erm
\end{bsp}

\begin{beme}
\brm
\begin{enum}
\item[]
\item In particular, a permutation $\pi\in{\cal S}_n(1243,2143)$ avoids $12\cdots k$ 
if and only if the path $\Psi_{EM}(\pi)$ has no step of height at least $k-1$. This result is equivalent to Theorem \ref{12...k-avoiding}.
\item Combining $\Psi_{EM}:{\cal S}_n(1243,2143)\to S_{n-1}$ with the bijection 
$\omega$ stated in Corollary \ref{bijection omega} yields the answer to the 
first part of the first problem raised in \cite{egge-mansour}. (By $S_{n-1}$ the 
set of Schr\"oder paths from $(0,0)$ to $(n-1,n-1)$ is denoted.) The map 
$\omega\circ\Psi_{EM}^{-1}:S_{n-1}\to{\cal S}_n(1243,2143)$ takes every 
Schr\"oder path whose maximum step height is at most $k-2$ to a $213\cdots 
k$-avoiding Schr\"oder permutation, and is bijective, of course.
\item Obviously, the path $\Psi_{EM}(\pi)$ contains no diagonal step if and 
only if ${\cal E}_0(\pi)=\emptyset$. As already noted, then $\pi_1=1$ and 
$\pi':=(\pi_2-1)(\pi_3-1)\cdots(\pi_n-1)$ belongs to ${\cal S}_{n-1}(132)$. In 
particular, we have $\Psi_{EM}(\pi)=\Psi_K(\pi')$ in this case.
\end{enum}
\erm
\end{beme}
\vspace*{0.5cm}


\setcounter{section}{5}\setcounter{satz}{0}

\centerline{\large{\bf 5}\hspace*{0.25cm}
{\sc Perspectives}}
\vspace*{0.5cm}

As already observed by Egge and Mansour in \cite{egge-mansour}, the 
investigation of $132$-avoiding permutations and $\{1243,2143\}$-avoiding ones, respectively, can be 
continued in a canonical way. For $m\ge3$ let $T_m$ be the set of permutations in 
${\cal S}_m$ for which $\pi_{m-1}=m$ and $\pi_m=m-1$. For example, 
$T_3=\{132\}$ and $T_4=\{1243,2143\}$. (The integer sequences counting the 
permutations in ${\cal S}_n(T_m)$ were determined in \cite{barcucci etal}.) Some of what we have done for 
$132$-avoiding permutations in \cite{reifegerste}, and for $T_4$-avoiding 
permutations in this paper can be generalized for an arbitrary integer $m$.

\begin{satz}
A permutation $\pi\in{\cal S}_n$ avoids each pattern in $T_m$ if and only if 
every element of its essential set is of rank at most $m-3$.
\end{satz}

\begin{bew}
If there exists an element $(i,j)\in{\cal E}(\pi)$ with $\rho(i,j)\ge m-2$ then at least 
$m-2$ dots appear northwest of 
$(i,j)$. Consequently, there are integers $i_1<i_2<\ldots<i_{m-2}<i$ for which 
$\pi_{i_1},\ldots,\pi_{i_{m-2}}<j$. Furthermore, we have $i<i_{m-1}$ and 
$\pi_i>j$ where $\pi_{i_{m-1}}=j$. Thus the subsequence 
$\pi_{i_1}\cdots\pi_{i_{m-2}}\pi_i\pi_{i_{m-1}}$ forms a pattern belonging to 
$T_m$. (For a better understanding draw a picture similar that one in the proof of Theorem \ref{characterization of Schroeder 
permutations}.)\\
On the other hand, it is clear from the diagram construction that the occurrence of a 
pattern of $T_m$ in a permutation yields a diagram corner of rank at least $m-2$.
\end{bew}

In order to study increasing subsequences in permutations which belong to ${\cal 
S}_n(T_m)$ it would be nice to have a surjection ${\cal S}_n(T_m)\to{\cal 
S}_n(T_{m-1})$ similar the map $\phi$ stated in Section 2. Thus the problem 
could successively put down to the case $m=3$.\\[2ex]
By reasoning similar to the proof of Theorem \ref{213...k-avoiding}, one can 
show that this theorem holds for each $m\ge3$ if $k=3$. In case $k\ge4$ and 
$m\ge5$, the condition $i+j\ge n+3-k+\rho(i,j)$ for all diagram corners $(i,j)$  
is only sufficient for avoiding $213\cdots k$ and all the patterns of $T_m$.\\ 
For example, the permutation $\pi=5\:4\:7\:1\:3\:2\:6\in{\cal S}_7(T_5)$ avoids 
$2134$ but the square $(1,4)$ is a diagram corner of rank $0$.\\[2ex]
Confirmed by computer tests, we believe in the Wilf-equivalence of $T_m\cup\{12\cdots k\}$ and 
$T_m\cup\{213\cdots k\}$ for all $k\ge1$ and $m\ge3$. 

\begin{verm}
For $m\ge 3$, and all $n$ and $k$ we have 
\bdpm
|{\cal S}_n(T_m\cup\{12\cdots k\})|=|{\cal S}_n(T_m\cup\{213\cdots k\})|.
\edpm
\end{verm}

\newpage


\centerline{\large\sc References}
\vspace*{0.5cm}

\begin{enumbib}

\bibitem{barcucci etal}
E. Barcucci, A. Del Lungo, E. Pergola, and R. Pinzani,  
Permutations avoiding an increasing number of length-increasing forbidden 
subsequences, 
{\it Discrete Math. Theor. Comput. Sci.} {\bf 4} (2000), 31-44. 

\bibitem{deutsch}
E. Deutsch, 
Dyck path enumeration, 
{\it Discrete Math.} {\bf 204} (1999), 167-202. 

\bibitem{egge-mansour}
E. S. Egge and T. Mansour, 
Permutations which avoid 1243 and 2143, continued fractions, and Chebyshev 
polynomials,
preprint, 2002.

\bibitem{eriksson-linusson}
K. Eriksson and S. Linusson, 
Combinatorics of Fulton's essential set,
{\it Duke Math. J.} {\bf 85} (1996), 61-80.

\bibitem{fulton}
W. Fulton, 
Flags, Schubert polynomials, degeneracy loci, and determinantal formulas,
{\it Duke Math. J.} {\bf 65} (1992), 381-420.

\bibitem{krattenthaler}
C. Krattenthaler, 
Permutations with restricted patterns and Dyck paths,
{\it Adv. Appl. Math.} {\bf 27} (2001), 510-530. 

\bibitem{kremer}
D. Kremer, 
Permutations with forbidden subsequences and a generalized Schr\"oder number, 
{\it Discrete Math.} {\bf 218} (2000), 121-130.

\bibitem{macdonald}
I. G. Macdonald, 
Notes on Schubert Polynomials,
LaCIM, Universit\'{e} du Qu\'{e}bec \`{a} Montr\'{e}al, 1991.

\bibitem{reifegerste}
A. Reifegerste, 
On the diagram of 132-avoiding permutations, 
preprint, 2002.

\bibitem{simion-schmidt}
R. Simion and F. W. Schmidt, 
Restricted Permutations, 
{\it Europ. J. Combinatorics} {\bf 6} (1985), 383-406. 

\bibitem{stanley}
R. P. Stanley, 
Enumerative Combinatorics Volume II,
Cambridge University Press, 1999.

\end{enumbib}
 
\end{document}